\documentclass[12pt]{amsart}

\newcommand{\R}{\mathbb{R}}
\newcommand{\N}{\mathbb{N}}

\newcommand{\Z}{\mathbb{Z}}

\newcommand{\T}{\mathbb{T}}

\newcommand{\A}{\mathcal{A}}

\newcommand{\MS}{M}

\renewcommand{\Pr}{\mathcal{P}}
\newcommand{\Qr}{\mathcal Q}
\newcommand{\Rr}{\mathcal R}

\newcommand{\E}{\mathcal{E}}
\newcommand{\hE}{\widehat{\E}}

\newcommand{\OP}{\Omega_{\Phi}}
\newcommand{\gmax}{\pi_{max}}
\newcommand{\Xmax}{X_{max}}
\newcommand{\Ymax}{X_{max}^{distal}}
\newcommand{\Phimax}{\Phi_{max}}

\newcommand{\larr}{\left( \begin{array}{c}}
\newcommand{\rarr}{\end{array} \right) }

\newcommand{\lsqarr}{\left[ \begin{array}{c}}
\newcommand{\rsqarr}{\end{array} \right]}

\newcommand{\inv}{\varprojlim}
\newcommand{\dir}{\varinjlim}

\newtheorem{theorem}{Theorem}[section]

\newtheorem{cor}[theorem]{Corollary}
\newtheorem{definition}[theorem]{Definition}
\newtheorem{corollary}[theorem]{Corollary}
\newtheorem{lemma}[theorem]{Lemma} \newtheorem{prop}[theorem]{Proposition}

\begin{document}

\title{Proximality and pure point spectrum for tiling
  dynamical systems}
\subjclass[2010]{Primary: 37B50, 37B05;
Secondary: 54H20, 54H11}
\keywords{Proximality, Maximal Equicontinuous Factor, Tiling Space}
\author{Marcy Barge} 
\address{Montana State University}
\email{barge@math.montana.edu}
\author{Johannes Kellendonk}
\address{Universit\'e de Lyon, Universit\'e Claude Bernard Lyon 1,
Institut Camille Jordan, CNRS UMR 5208, 43 boulevard du 11 novembre
1918, F-69622 Villeurbanne cedex, France}
\email{kellendonk@math.univ-lyon1.fr}

\date{\today}

\begin{abstract} 
We investigate the role of the proximality relation for tiling
dynamical systems. Under two hypothesis, namely that the minimal rank
is finite and the set of fiber distal points has full measure we show that
the following three conditions are equivalent: (i) proximality is
topologically closed, (ii) the minimal rank is one, (iii) the
continuous eigenfunctions of the translation action span the $L^2$-functions
over the tiling space.  We apply our findings to model sets and to Meyer substitution tilings. 
It turns out that
the Meyer property is crucial for our analysis as it allows us to
replace proximality by the a priori stronger notion of strong proximality.
\end{abstract} 

\maketitle

\tableofcontents

\section{Introduction}\label{intro} In order to understand the
combinatorial properties of a single tiling $T$ of Euclidean space
$\R^n$, one may consider the collection $\Omega_T$, 
called the hull of $T$,
of all tilings of $\R^n$ that are locally indistinguishable from $T$:
a tiling $T'$ is in the hull of $T$ if every finite collection of
tiles in $T'$ is exactly a translate of a finite collection of tiles
of $T$. There is a natural topology on the hull, and an action of
$\R^n$ on the hull by translation, so that the properties of the
resulting topological dynamical system reflect combinatorial
properties of the original tiling. A beautiful example of this
correspondence between combinatorics and dynamics arises in
diffraction theory. The diffraction spectrum of a point set in $\R^n$
(think of the points as atoms in a material and the
diffraction spectrum as a picture of X-ray scattering) depends on the
spatial recurrence properties of finite patterns of points in the
point set. The point set determines a tiling of $\R^n$ by Veronoi
cells and, if the patches of the tiling are distributed in a sufficiently  
regular manner, 
the point set has pure
point diffraction spectrum (the material is a perfect quasicrystal) if
and only if the $\R^n$-action on the hull of the tiling has pure
discrete dynamical spectrum (\cite{D}, \cite{LMS}). The latter means
by definition 
that the Hilbert space $L^2(\Omega_T,\mu)$ is generated by the
continuous eigenfunctions of the action\footnote
{Eigenfunctions are (classes of) functions
$f:\Omega_T\to \mathbb C$ satisfying $f(T'-v) = \exp(2\pi\imath \beta(v))f(T')$
for all $v\in\R^n$ almost all $T'\in\Omega_T$ and some linear
functional $\beta:\R^n\to\R$ (the eigenvalue).}. 
Here $\mu$ is an invariant
ergodic measure on $\Omega_T$ which is related to the diffraction via the 
construction of the autocorrelation measure\footnote{We will consider
below strictly 
ergodic tiling dynamical in which case this measure is unique.}.

The study of continuous eigenfunctions is related to the study of
equicontinuous factors of the dynamical system $(\Omega_T,\R^n)$. 
All
continuous eigenfunctions together determine what is called the
maximal equicontinuous factor $\gmax:(\Omega_T,\R^n)\to (\Xmax,\R^n)$.
One of the most common routes to determine whether  
the $\R^n$-action on the hull of the tiling has pure
discrete dynamical spectrum is therefore 
to examine whether $\gmax$ is almost everywhere one-to-one. 

The equivalence relation whose equivalence classes 
are the fibers of $\gmax$
is called the equicontinuous structure relation. It is related
to the proximal relation. The latter is not necessarily an equivalence
relation but Auslander (\cite{Aus}) has shown that a modification of proximality 
gives rise to a relation which he called regional proximality and
which, for minimal systems, coincides with the equicontinuous
structure relation. Our main aim in this work is to use these concepts of proximality to
say something about the 
equicontinuous structure relation, in particular in light of the
above question whether the dynamical spectrum of a tiling system is pure point.

The concept of proximality applies to general topological dynamical
systems $(X,G)$ and the definition given in Auslander's book 
requires only that $X$ be a compact uniform space and
carry a continuous group $G$ action $\alpha$. Since the hull of a tiling is
metrizable we can work with a metric $d$. This involves an irrelevant
choice of metric but is a little less abstract. So we
consider a compact metric space $(X,d)$ with a minimal 
continuous $G$ action. We require $G$ to be a locally compact abelian
group and for some results also that $G$ is compactly generated.

Two points $x,y\in X$ are proximal if
$\inf_{t\in G} d(\alpha_t (x),\alpha_t( y)) = 0$. A point $x\in X$ is distal
if it is not proximal to any other point.
We say that $x\in \Xmax$ is fiber distal if $\gmax^{-1}(x)$ consists only of distal points.
Let $\gmax:(X,G)\to (\Xmax,G)$ be the maximal equicontinuous factor. 
We say that $(X,G)$ has finite minimal rank if its minimal rank
$$ mr := \inf\{\#\gmax^{-1}(x):x\in \Xmax\}$$
is finite.
We show
\medskip

\noindent
{\bf Theorem} (Theorem~\ref{thm-P=Q} in the main text)
{\em Let $(X,G)$ have finite minimal rank and suppose that $G$ is compactly generated.
Then the proximal relation $\Pr$
coincides with the equicontinuous structure relation if and only if 
$\Pr\subset X\times X$ is closed (in the product topology).}
\medskip

\noindent
{\bf Theorem} (Lemma~\ref{lem-cr-mr} together with  Lemma~\ref{lem-cr=mr} in the main text)
{\em Let $(X,G)$ have finite minimal rank and suppose that its maximal equicontinuous factor
admits at least one fiber distal point. 
Then the proximal relation $\Pr$ is closed
if and only if the minimal rank is $1$.}
\medskip

Now let $\eta$ be the (normalised) Haar measure on $\Xmax$. 
We say that $(X,G)$ is almost everywhere fiber distal if the set of
fiber distal points has full measure in $\Xmax$. 
\medskip

\noindent
{\bf Theorem} (Theorem~\ref{thm-1.3} in the main text)
{\em Let $(X,G)$ have finite minimal rank and be almost everywhere fiber distal.Let $\mu$ be an ergodic $G$-invariant  Borel probability measure on $X$. 
Then $L^2(X,\mu)$
is generated by continuous eigenfunctions  
if and only if the proximal relation is closed.}
\medskip

In particular, if we know already that all eigenfunctions are
continuous then, under the hypothesis that $(X,G)$ has finite minimal
rank and the fiber distal points have full measure, topological closedness of $\Pr$ is
a neccessary and sufficient criterion for pure discrete spectrum.
It is therefore of interest to investigate which type of tiling (or Delone
set) systems satisfy the hypotheses.
We show here that this is the case for regular model sets and Meyer substitution 
tilings.
\bigskip

The hull of a tiling has a laminated structure which comes from
the group action. Unlike more general dynamical systems, this lamination admits canonical
transversals. The extra transverse structure allows for the definition of a stronger
notion of proximality: we call two tilings strongly proximal if they
agree exactly on arbitrarily large balls. It is central to our analysis that
for the most important class of Meyer sets both notions coincide. A
similar result applies to the regional proximal relation: there is a
strong version of it which coincides with the usual one for Meyer
sets. 

A large part of the paper is devoted to the study of two classes of Meyer systems, those defined by model
sets and those defined by Meyer substitutions. We find
\begin{itemize}
\item Model sets always have minimal rank $1$ and proximality is thus always closed. The set of fiber distal points has full measure if and only if the model set is regular.
\item Meyer substitution tilings always have finite minimal
  rank and the set of fiber distal points always has full measure. 
\end{itemize}
The case of model sets seems a lot simpler, but, except for nice windows, we cannot control
the maximal rank, $\sup\{\#\gmax^{-1}(x):x\in \Xmax\}$. On the
other hand, we find that Meyer substitution tilings always have
finite maximal rank. This is tremendously advantageous also in other
contexts.  It is exploited in \cite{BO} to describe
the branch locus in 2-$d$ self-similar Pisot substitution tiling
spaces  
and in \cite{B} to characterize minimal directions in self-similar
Pisot substitution tiling spaces of any dimension.

We consider a third notion of proximality. Two points $x,y$ are completely proximal if, given any subset $A\subset G$ which contains a translate of each compact subset, we have
$\inf_{t\in A} d(\alpha_t(x),\alpha_t(y)) = 0$.
This turns out to be an equivalence
relation, but we do not know whether it is always closed. We show that for Meyer substitutions complete proximality is indeed a closed equivalence relation.
All three notions - complete
proximality, proximality and regional proximality - are the same for
tilings of finite minimal rank if proximality is closed. Thus for Meyer substitutions we obtain 
(Cor.~\ref{cor-Meyer} in the main text) equivalence between:
\begin{enumerate}
\item proximality is a closed relation, 
\item proximality agrees with complete proximality,
\item the dynamical spectrum is purely discrete.
\end{enumerate}
\section{Maximal equicontinuous factors and
  proximality}\label{Maximal}

\subsection{General notions and results}
In this section we recall some aspects of the theory of topological dynamical systems relating to equicontinuity, proximality and regional proximality. This is mainly based on the material of Auslander's book
\cite{Aus}.   

We consider a dynamical system $(X,G)$ where $X$ is a compact
metrizable space and $G$ a locally compact abelian group acting
continuously by $\alpha$ on $X$. We denote the action by
$\alpha_t(x) = t\cdot x$, or, in the context of tilings, by $\alpha_v(T)=T-v$. 
\begin{definition}[equicontinuity]  A point $x\in X$ is
  called equicontinuous if the family of homeomorphisms $\{\alpha_t \}_{t\in G}$ is
  equicontinuous at $x$. The dynamical system $(X,G)$ is called equicontinuous
  if all its points are equicontinuous.
\end{definition}
Although the standard definition of equicontinuity uses a metric, it does not depend on the particular choice of metric, as long as the metric is compatible with the topology. In fact, Auslander introduces this notion using the (unique) uniformity defined by the topology.
Minimal equicontinuous systems have a very simple structure: they are translations on compact abelian groups. This means that $X$ has the structure of an abelian group and the action is given by 
$\alpha_t (x) = x+\imath(t)$ where $\imath:G\to X$ is a group homomorphism.
\begin{theorem}[Ellis]\label{thm-equi}
A minimal system $(X,G)$  is equicontinuous if and only if it
is conjugate to a minimal translation on a compact abelian group.
\end{theorem}
If $(X,G)$ is equicontinuous the group structure on $X$ arises as follows: given any point $x_0\in X$ the operation $t_1\cdot x_0 + t_2\cdot x_0 := (t_1 + t_2)\cdot x_0$ extends to an addition in $X$ so that
$X$ becomes a group with $x_0$ as neutral element. Conversely any translation on a compact abelian group is clearly equicontinuous.

We are interested in dynamical systems defined by aperiodic tilings of finite local complexity (FLC). Such systems are never equicontinuous \cite{BO} and it will prove fruitful to study the relation between these systems and their maximal equicontinuous factors. 
\begin{definition} [maximal equicontinuous factor]  An equicontinuous factor of $(X,G)$ is maximal if any other equicontinuous factor of $(X,G)$ factors through it. It is thus unique up to conjugacy and therefore referred to as {\bf the} maximal equicontinuous factor. We denote it $(\Xmax,G)$ and the factor map by $\gmax:X\to \Xmax$.\end{definition}
The maximal continuous factor always exists but may be trivial ({\it i.e.}, a single point).  The equivalence relation defined by $\gmax$, that is, $x\sim y$ if $\gmax(x)=\gmax(y)$, is called the {\bf equicontinuous structure relation}.

The concept of proximality is central to Auslander's investigation of equicontinuous structure relation.
\begin{definition}
[proximality] Consider a compatible metric $d$ on $(X,G)$.
Two points $x,y\in X$ are proximal if $$\inf_{t\in G} d(t\cdot x,t\cdot y)=0.$$
We denote by $\Pr\subset X\times X$ the proximal relation and write $x\sim_py$ if $(x,y)\in \Pr$.
\end{definition}
The proximal relation does not depend on the metric but only on the topology (it can as well be formulated using the uniformity structure on $X$). It is easy to see that the proximal relation is trivial for equicontinuous systems, but
the converse is not true. Systems for which the proximal relation is trivial are called {\bf distal}.

Note that $\Pr = \bigcap_\epsilon \Pr_\epsilon$ where 
$$\Pr_\epsilon = \{(x,y) \in X\times X: \inf_{t\in G} d(t\cdot x,t\cdot y)<\epsilon\}.$$
The proximal relation is not always closed, i.e.,\ $\Pr$ need not be a closed subset of $X\times X$,
and is not, in general, a transitive relation. 
However:
\begin{theorem}[{\cite{Aus}}]{\label{thm-Aus1}}
If the proximal relation is closed then it is an equivalence relation.
\end{theorem}

\begin{definition}
[regional proximality] 
The regional proximal relation is $\Qr:=\bigcap_{\epsilon}\overline{\Pr_\epsilon}$. In other words,
two points $x,y\in X$ are regionally proximal if for all $\epsilon$ there exist $x'\in X$, $y'\in X$, $t\in G$ such that
$d(x,x')<\epsilon$, $d(y,y')<\epsilon$ and $d(t\cdot x',t\cdot y')<\epsilon$.
\end{definition}
\begin{theorem}[{\cite{Aus}}] The equicontinuous structure relation is the smallest closed equivalence relation
containing the regional proximal relation.
\end{theorem}
The regional proximal relation is, in general, neither closed nor transitive. 
However, if the acting group is abelian, we have:
\begin{theorem}[{\cite{Aus}}] For minimal systems the regional proximal relation is 
a closed equivalence relation and hence it coincides with the equicontinuous structure relation.
\end{theorem}
Even for minimal systems, the regional proximal relation is not necessarily the smallest closed equivalence
relation containing the proximal relation.
 Indeed, $\Pr$ may be trivial while $\Qr$ is not (there are minimal distal systems which are not equicontinuous). 
 \bigskip
 
The following definition is a generalisation of a notion which has been introduced in the context of Pisot substitution tilings in \cite{BK}.
For $\delta>0$ and $x\in \Xmax$ let 
$ cr(x,\delta)$ be  the maximal cardinality $l$ of a collection
$\{x_1,\ldots,x_l\}\subset \gmax^{-1}(x)$  
with the property that $\inf_{t\in G} d(t\cdot x_i,t\cdot x_j)\geq\delta$
provided $i\neq j$. 

\begin{definition}[coincidence rank] The coincidence rank of a minimal system  $(X,G)$ is the number
$$cr = \lim_{\delta\to 0^+} cr(x,\delta).$$ 
\end{definition}
\begin{lemma} \label{lem-14} 
The limit in the above definition does not depend on the choice of $x$.
Moreover, if $cr(x,\delta)$ is finite for some $x$, then 
there exists a $\delta_0>0$ such that, for all $y$, $$\lim_{\delta\to 0^+} cr(y,\delta)= cr(y,\delta_0).$$ 
\end{lemma}
\begin{proof}

Clearly, $cr(x,\delta)$ is a decreasing integer-valued function of
$\delta$ and $cr(x,\delta) = 1$ if $\delta$ is larger than the
diameter of $X$. Furthermore, $cr(x,\delta) =cr(t\cdot x,\delta) $ for all
$t\in G$.  
This implies that either $\lim_{\delta\to 0^+}cr(x,\delta)= +\infty$ or 
$\lim_{\delta\to 0^+}cr(x,\delta)= cr(x_0,\delta_0)$ for some
$\delta_0>0$ and all points $x$ of the orbit of  a point $x_0\in
\Xmax$.  
We need to show that the result is the same for points $y\in \Xmax$ in
other orbits. 

Let $l\in\N$ and $\{x_1,\ldots,x_l\}\subset \gmax^{-1}(x_0)$ 
with the property that $\inf_{t\in G} d(t\cdot x_i,t\cdot x_j)\geq\delta$
provided $i\neq j$.  
Let $x'_1\in \gmax^{-1}(y)$. By transitivity, there exists a sequence
$(t_n)_n$ such that $\lim_n t_n\cdot x_n\to x'_1$. Taking subsequences we
may suppose that all other limits  $\lim_n t_n\cdot x_i$ exist: let $x'_i$
denote these limits. 
Note that since $\gmax$ is continuous and equivariant w.r.t.\ the
action we have $x'_i\in\gmax^{-1}(y)$ for all $i=1,\cdots,l$. 
We have
$$ d(t'\cdot x'_1,t'\cdot x'_2) \geq  d(t'\cdot (t_n\cdot x_1),t'\cdot (t_n\cdot x_2)) -
d(t'\cdot x'_1,t'\cdot (t_n\cdot x_1)))- d(t'\cdot (t_n\cdot x_2),t'\cdot x_2).$$ 
Keeping $t'$ fixed we find, for all $\epsilon>0$, an $n$ such that 
$d(t'\cdot x'_1,t'\cdot (t_n\cdot x_1))<\epsilon$ and $d(t'\cdot (t_n\cdot x_2),t'\cdot x'_2)<\epsilon$. Hence
$ d(t'\cdot x'_1,t'\cdot x'_2) \geq \delta - 2\epsilon$ and, since $\epsilon$ was
arbitrary, we see that  
$cr(y,\delta)\geq cr(x,\delta)$. By a symmetric argument,
$cr(x,\delta)\geq cr(y,\delta)$.  
\end{proof}
\begin{definition}[minimal rank] The
minimal rank of $(X,G)$ is $$mr:= \inf\{\#\gmax^{-1}(x):x\in\Xmax\}.$$
\end{definition}
\begin{lemma}\label{lem-cr-mr} 
$cr=1$ if and only if $\Pr = \Qr$. Furthermore $cr\le mr$. 
\end{lemma}
\begin{proof} 
Suppose that $\Pr = \Qr$. Then for all $x_1,x_2\in \gmax^{-1}(x)$ we have
$\inf_{t\in G} d(t\cdot x_1,t\cdot x_2)=0$ and hence $cr(x,\delta) = 1$. 
 
Now suppose that $cr=1$. Since $cr(x,\delta)$ is a decreasing function
of $\delta$ we have 
$cr(x,\delta) = 1$ for all $x$ and $\delta>0$. In particular, two
elements $x_1,x_2 \in\gmax^{-1}(x)$ cannot satisfy 
$\exists \delta>0:\inf_{t\in G} d(t\cdot x_1,t\cdot x_2)\geq \delta$. The
latter means that $(x_1,x_2)\notin \Pr_\delta^c$ 
for all $\delta>0$. In other words $(x_1,x_2)\in \bigcap_\delta \Pr_\delta$. 

The inequality is clear from the independence of $\lim_{\delta\to 0^+}
cr(x,\delta)$ of $x$. 
\end{proof}
\begin{definition}[distal, fiber distal]
A point $x\in X$ is distal if for all $y\in X$, 
$\inf_{t\in G} d(t\cdot x,t\cdot y)>0$, i.e.,\ $x$ is not proximal to any other
point. We say that a point $x\in\Xmax$ is fiber distal, if all points
of the fiber $\gmax^{-1}(x)$ are distal.
We denote by  $ \Ymax\subset \Xmax$ the set of fiber distal points. 
$(X,G)$ is called fiber distal if its maximal equicontinuous factor admits
 a fiber distal point.
\end{definition}
Note that, since $\Pr$ is contained in the equicontinuous structure
relation a point $x$ is distal if and only if it is 
not proximal to any other point in the fiber $\gmax^{-1}(\gmax(x))$. 
\begin{lemma}\label{lem-cr=mr} 
Let $(X,G)$ be a minimal system with finite minimal rank.
Then $cr=\#\gmax^{-1}(x)$ whenever $x$ is fiber distal.
In particular, $ \Ymax = \{x\in\Xmax:\#\gmax^{-1}(x)=cr\}$ and
 $cr=mr$ whenever $\Xmax$ contains a fiber distal point. 
\end{lemma}
\begin{proof}
Let $x\in \Xmax$ be a fiber distal point and $\{x_1,\cdots,x_k\} \subset
\gmax^{-1}(x)$, $k\geq mr$. 
Then none of the points $x_i$ is proximal to
any other $x_j$. By definition this means that
$cr(\gmax(x),\delta_0)\geq k$ for $0<\delta_0 <\inf_{i\ne j\in\{1,\ldots,k\}}  \inf_{t\in G}
d(t\cdot x_i,t\cdot x_j)$. Hence $cr\geq k \geq mr$. We have already seen that $mr\geq
cr$. This shows also that $cr=\#\gmax^{-1}(x)$ if $x$ is fiber distal.

Now suppose that $cr=mr$. Then there exists $\xi\in \Xmax$ such that
$\gmax^{-1}(\xi) = \{x_1,\cdots,x_{mr}\}$. Since $cr =
cr(\xi,\delta_0)$ for some $\delta_0$, we must have that  $\inf_{t\in G}
d(t\cdot x_i,t\cdot x_j)\geq \delta_0$ for all $i\neq j$. Thus all $x_i$ are distal.

Finally, if  $cr=\#\gmax^{-1}(x)$ then $\gmax^{-1}(x)$
cannot contain proximal points and so
$x$ must be fiber distal. 
\end{proof}

Let $\tilde \Pr$ the smallest closed equivalence relation containing the
proximal relation on a minimal compact metrizable dynamical system
$(X,G)$. We define $X_p$ to be the quotient space $X_p:= X/\tilde \Pr$. Since
$\tilde \Pr$ is closed, 
$X_p$ is metrizable and the canonical projection is a closed continuous
map \cite{Kurka}.  
Since the proximal relation is $G$-invariant the action of $G$ descends and so we have a factor system 
$(X_p,G)$.  Furthermore, $\gmax$ factors through the above canonical projection and so we get another factor map    
$$ \pi : X_p \to \Xmax,\quad  \pi([x]_p) = \gmax(x),$$
which is again closed.

\begin{theorem}\label{thm-P=Q} Let $(X,G)$ be a minimal system with finite minimal rank. Suppose that $G$ is compactly generated. Then $\Pr=\Qr$ if and only if $\Pr$ is closed.
\end{theorem}

The proof of the theorem is based on the following two lemmas. We know already that $\Pr=\Qr$ is equivalent to $cr=1$. What we need to show, therefore, is that $\Pr$ closed implies $cr=1$ (under the assumption that $cr$ is finite!).

\begin{lemma}
Suppose that $cr$ is finite. If $\Pr$ is closed, then $\pi$ is a $cr$-to-$1$ map which is a local homeomorphism, i.e.,\
any point in $X_p$ admits a neighborhood on which $\pi$ restricts to homeomorphism onto its image.
\end{lemma}
\begin{proof}
Note that if $\Pr$ is closed, then $\Pr=\tilde{\Pr}$ by Theorem \ref{thm-Aus1}. We first show that $\pi$ is $cr$-to-$1$. For $x\in\Xmax$, clearly $\#\pi^{-1}(x) \geq cr$ as there are elements $x_1,\cdots,x_{cr}\in \gmax^{-1}(x)$ that belong to different $\Pr$-classes. Suppose now that  $x_1,\cdots,x_{l}\in \gmax^{-1}(x)$ with $x_i\not\sim_p x_j$ for all $x_i\neq x_j$. Then there exists $\delta>0$ such that  for all $x_i\neq x_j$,
$(x_i,x_j)\notin \Pr_\delta$. By definition, $l\leq c(x,\delta)\leq cr$.

Clearly $\pi$ is continuous and surjective.
If $\pi$ were not locally injective, we could construct two sequences $(\xi_n)_n,(\eta_n)_n$ in $X/\sim_p$ such that,
\begin{enumerate}
\item\label{item-1} $\xi_n\neq \eta_n$ for all $n$,
\item $\pi(\xi_n) = \pi(\eta_n)$ for all $n$,
\item $\lim_n \xi_n = \lim_n\eta_n = [x]_p$, for some $x\in X$.
\end{enumerate}
Suppose there are such sequences  $\xi_n = [x_n]_p$, $\eta_n=[x'_n]_p$. We may suppose without loss of generality that $\lim_n x_n = x$ and $\lim x'_n = x'$ for some $x'\sim_p x$. Then $\inf_{t\in G} d(t\cdot x,t\cdot x') = 0$.
For all $\epsilon >0$ there exists $N_{\epsilon,t}$ such that $\forall n\geq N_{\epsilon,t}$ we have
$d(t\cdot x_n, t\cdot x)<\epsilon$ and $d(t\cdot x'_n, t\cdot x')<\epsilon$. So if $t$ is such that $d(t\cdot x,t\cdot x') < \epsilon$ and 
$n\geq N_{\epsilon,t}$, we have $d(t\cdot x_n,t\cdot x'_n)<3\epsilon$. Now if $3\epsilon<\delta_0$ (from Lemma \ref{lem-14}) then $x_n\sim_p x'_n$ for all $n\geq N_{\epsilon,t}$. This violates (\ref{item-1}).

We thus have shown that $\pi$ is a locally injective continuous surjection. Furthermore $\pi$ is a closed map,
hence the restriction of $\pi$ to an open neighborhood, on which $\pi$ is injective, yields a closed continuous invertible map onto the image of that neighborhood under $\pi$. The restriction is therefore a homeomorphism.  \end{proof}

\begin{lemma} Let $(Y,G)$ be a compact metrizable dynamical system and
$\pi:(Y,G)\to (X,G)$ a finite-to-one factor map which is a local homeomorphism.
We suppose that $X$ is connected and that
$G$ is a compactly generated abelian group, i.e.,\
there is a compact neighborhood $K$ of the identity $0$ in $G$ so that $G=\cup_{n\in\N}nK$, $nK:=K+\cdots+K$.
If $(X,G)$ is equicontinuous then $(Y,G)$ is also equicontinuous.
\end{lemma}
\begin{proof}
Since $X$ is connected and $Y$ is compact, there is $m\in\N$ so that $\pi$ is $m$-to-$1$ everywhere. In the following, $\{W_i\}$ will denote a finite open cover of $X$ with the property that $\pi^{-1}(W_i)$
is the disjoint union of $W_{i,j}$, $j=1,\ldots,m$, with $\pi|_{W_{i,j}}:W_{i,j}\to W_i$ a homeomorphism for all $i,j$. We may take the $W_i$ small enough so that $\pi|_{\bar{W}_{i,j}}:\bar{W}_{i,j}\to \bar{W}_i$ is a homeomorphism for each $i,j$, with inverse denoted by $\pi_{i,j}^{-1}$. Let $\alpha>0$ be a Lebesgue number for both $\{W_i\}$ and $\{W_{i,j}\}$. 
We denote by $d_X$ and $d_Y$ chosen metrics on $X$ and $Y$. Since $(X,G)$ is equicontinuous we may suppose that $d_X$ is invariant under the action of $G$. Let 
 $$\beta:=\min_{i, j\neq k} d_Y(\bar W_{i,j},\bar W_{i,k}) = \min\{d_Y(y_1,y_2):y_1\in\bar{W}_{i,j},\, y_2\in\bar{W}_{i,l},\, l\ne j\}.$$
 Since $K$ is compact there is 
$\gamma>0$ so that if $d_Y(y_1,y_2)<\gamma$ then $d_Y(t\cdot y_1,t\cdot y_2)<\beta$  for all $t\in K$.

Now let $0<\epsilon\leq \gamma$ be given. 
\begin{itemize}
\item Let $0<\eta\leq \alpha$ so that if $d_X(x_1,x_2)<\eta$ and $x_1,x_2\in W_i$, then  for all $j$
$d_Y(\pi_{i,j}^{-1}(x_1),\pi_{i,j}^{-1}(x_2))<\epsilon$. 
\item Let $0<\delta<\min\{\alpha,\gamma\}$, so that if $d_Y(y_1,y_2)<\delta$, then $d_X(\pi(y_1),\pi(y_2))<\eta$. 
\end{itemize}
We prove with an inductive argument that if $d_Y(y_1,y_2)<\delta$, then $d_Y(t\cdot y_1,t\cdot y_2)<\epsilon$ for all $t\in G$. 

First suppose that $t\in K$ and $d_Y(y_1,y_2)<\delta$. Set $x_1=\pi(y_1)$ and $x_2=\pi(y_1)$. 
Then $d_X(t\cdot x_1,t\cdot x_2)=d_X(x_1,x_2)< \eta\leq \alpha$ so there there is $i$ with $t\cdot x_1,t\cdot x_2\in W_i$ and there are $j,l$ with $t\cdot y_1\in W_{i,j}$ and $t\cdot y_2\in W_{i,l}$.
Since $d_Y(y_1,y_2)<\gamma$ and $t\in K$, it must be that $j=l$. Thus $t\cdot y_1=\pi_{i,j}^{-1}(t\cdot x_1)$
and $t\cdot y_2=\pi_{i,j}^{-1}(t\cdot x_2)$. Since $d_X(t\cdot x_1,t\cdot x_2)<\eta$ we have $d_Y(t\cdot y_1,t\cdot y_2)<\epsilon$. To summarize, we have shown that $t\in K$ and $d_Y(y_1,y_2)<\delta$ implies 
$d_Y(t\cdot y_1,t\cdot y_2)<\epsilon$.

Now assume, for some $n\in\N$, that
$d_Y(t\cdot y_1,t\cdot y_2)<\epsilon\leq\gamma$ for all $t\in nK$ and all
$y_1,y_2$ with $d_Y(y_1,y_2)<\delta$. Let $t=t_1+t_2$ with $t_1\in nK$
and $t_2\in K$, and suppose that $d_Y(y_1,y_2)<\delta$. Let
$y_1':=t_1\cdot y_1$ and $y_2':=t_1\cdot y_2$. Then
$d_Y(y_1',y_2')<\gamma$ and $d_X(t\cdot x_1,t\cdot x_2)=d_X(x_1,x_2)<\eta$. As
above, we conclude that there are $i,j,l$ with $t_2\cdot y_1'\in W_{i,j}$ and
$t_2\cdot y_2'\in W_{i,l}$. Since $t_2\in K$, $j=l$. Then
$t_2\cdot y_1'=\pi_{i,j}^{-1}(t\cdot x_1)$ and  
$t_2\cdot y_2'=\pi_{i,j}^{-1}(t\cdot x_2)$, so that $d(t\cdot y_1,t\cdot y_2)<\epsilon$.

By induction, if $d(y_1,y_2)<\delta$, then $d(t\cdot y_1,t\cdot y_2)<\epsilon$ for all $t\in\cup_{n\in\N}nK=G$.
\end{proof}

\begin{proof}(of Theorem~\ref{thm-P=Q})
Combining the last two lemmas we find that $(X_p,G)$ is equicontinuous and hence must already coincide with 
$(\Xmax,G)$. Thus $\pi$ is a conjugacy and $cr=1$.
\end{proof}
The above is all we will need for our analysis of tiling systems. For completeness, we add a corollary about distal systems and strengthen the last theorem. A dynamical system is called distal if all its points are distal.
\begin{cor}
For compact metrizable minimal distal systems we have either $cr=1$ or $cr=+\infty$,
the first case arising if and only if the system is equicontinuous.
\end{cor}
Note that the following statements are equivalent for a monoid $E$
(a semi-group with neutral element):
\begin{enumerate}
\item $E$ is a group,
\item any element of $E$ admits a left inverse,
\item $E$ does not admit non-trivial ideals. 
\end{enumerate}
\begin{prop} Let $(X,G)$ be a compact Hausdorff dynamical system.
$(X_p,G)$ is distal.
\end{prop}
\begin{proof} We need to  prove that the Ellis monoid $E(X_p)$ is a group.
The canonical projection $\pi:X\to X_p$ induces a continuous semigroup homomorphism
$\pi_*:E(X)\to E(X_p)$ determined by the equation
$$ \pi_*(p)(y) = \pi(p(x))$$
where $x$ is any preimage of $y$. Suppose that $E(X_p)$ is not a group and therefore admits a non-trivial ideal $I$.
Then $\pi_*^{-1}(I)$ is a non-trivial closed ideal of $E(X)$ (closed by continuity). Auslander (\cite{Aus}) shows that any closed ideal of the Ellis semigroup of a dynamical system contains an idempotent $u$. Any idempotent $u\in E(X)$ satisfies $u(x) \sim_p x$ for all $x\in X$. Hence $\pi(u(x)) = \pi(x)$, so that $\pi_*(u)(x) = \pi(u(x)) = \pi(x)$ for all $x\in X$. Thus $\pi_*(u) = \mbox{id}$ and $\mbox{id}\in I$, contradicting that $I$ is proper.  
\end{proof}
\begin{cor}
Consider a minimal system $(X,G)$ of finite minimal rank. Then the equicontinuous structure relation is the smallest closed equivalence relation containing the proximal relation.
\end{cor}
\begin{proof}
Since the minimal rank of $(X,G)$ is finite and $(X_p,G)$ is a factor sitting above the equicontinuous factor
the minimal rank of $(X_p,G)$ must be finite too. By the last result this implies that $(X_p,G)$ is equicontinuous and hence coincides with the maximal equicontinuous factor.
\end{proof}


\subsection{The maximal equicontinuous factor and the dynamical spectrum}
We have seen above that the maximal equicontinuous factor arises from
dividing out the regional proximal relation (or the closed equivalence
relation generated by it, in the non minimal case). 
There is another way to describe the maximal equicontinuous
factor and this is related to the topological dynamical spectrum of the system.

A {\bf continuous eigenfunction} of a dynamical system $(X,G)$ is a non-zero function
$f\in C(X)$ for which there exists a (continuous) character $\chi\in
\hat G$ such that  
\begin{equation}\label{eq-eigen}
f(t\cdot x) = \chi(t) f(x).
\end{equation} 
$\chi$ is called the {\bf eigenvalue} of $f$ and the set of all eigenvalues
$\mathcal E$ forms a subgroup of the Pontrayagin dual $\hat G$ of $G$.  
We call the eigenfunction normalized if its modulus is equal to $1$.

Note that, by universality of the maximal equicontinuous factor, all continuous eigenfunctions on $(X,G)$ factor through $\gmax$. Indeed, any normalised continuous eigenfunction $f:X\to \T^1$ gives rise to an equicontinuous factor of the form $(\T^1,G)$ which hence sits below the maximal equicontinuous factor, i.e.\ $f=f'\circ \gmax$ for some $f'\in C(\Xmax)$. In particular, $\gmax(x)=\gmax(y)$ implies that all  continuous eigenvalues take the same value on $x$ as on $y$. 
Given that the maximal equicontinuous factor is a minimal translation on a compact abelian group, its eigenvalues are $\widehat\Xmax$ and the continuous eigenfunctions separate the points of $\Xmax$. Now If $\gmax(x)\neq \gmax(y)$ then there exists an continuous eigenfunction $\tilde f$ for $(\Xmax,G)$ such that $\tilde f(\gmax(x))\neq \tilde f(\gmax(y))$. Hence $\gmax^*(\tilde f)$ is a continuous eigenfunction taking different values on $x$ and $y$.  Combining the two arguments we see that
$\gmax(x)=\gmax(y)$ if and only if for all continuous eigenfunctions $f$ one has $f(x)=f(y)$.

Let $\mathcal F$ be the norm closed sub-algebra of
$C(X)$ generated by the continuous eigenfunctions. We will argue that
the maximal equicontinuous factor $\Xmax$ can be identified with
the spectrum\footnote{The spectrum is  the space of non-zero
$*$-algebra morphisms $\varphi:\mathcal F\to \mathbb C$ equipped with
the subspace topology of the weak-*-topology of the dual space
$\mathcal F^*$.}
$\hat{\mathcal F}$ of $\mathcal F$. The action of $G$ on $C(X)$ via pull-back
preserves the space of continuous eigenfunctions and hence gives
rise by duality to an action on $\hat{\mathcal F}$: for an eigenfunction $f$ with eigenvalue $\chi$ this is
$$(t\cdot \varphi)(f): = \chi(t) \varphi(f).$$
The dual of
the inclusion $\imath :\mathcal F\to C(X)$ yields therefore a factor
map $\hat\imath:\widehat{C(X)}\cong X\to \hat{\mathcal F}$. Taking into account the homeomorphism $x\mapsto ev_x$ between $X$ and the spectrum $\widehat{C(X)}$ of $C(X)$ (the latter is given by the set of evaluation maps $\{ev_x:x\in X\}$, $ev_x(f):=f(x)$) the map $\hat\imath(x)$ is simply the 
restriction of $ev_x$ to the space of continuous eigenfunctions.
By the above remarks
$\gmax(x)=\gmax(y)$ if and only if $\hat\imath(x) = \hat\imath(y)$ and hence
$\gmax$ and $\hat\imath$ have the same fibers. It follows that 
the factor $\hat\imath:X\to \hat{\mathcal F}$ is
isomorphic to the maximal equicontinuous factor.

Finally, observe that $\hat{\mathcal F}$ can be identified with the
Pontrayagin dual $\hat{\mathcal E}$ of $\mathcal E$. Indeed, 
choose a point $x_0\in X$ and normalise all continuous eigenfunctions to
$f(x_0)=1$. Then there is a bijection between normalized
continuous eigenfunctions and eigenvalues\footnote{We assume minimality here,
  otherwise one has to choose several points.}; 
namely, $\chi\leftrightarrow f_\chi$, since
$f_\chi (t\cdot x_0) :=\chi(t)$ extends by continuity to the unique normalized continuous eigenfunction to eigenvalue $\chi$. 
This is even a abelian group
isomorphism. It follows that, if we equip $\mathcal E$ with the
discrete topology, then the continuous eigenfunctions form the group algebra $\mathbb C \mathcal E$ 
and hence $\hat{\mathcal F}=\hat{\mathcal E}$, a compact abelian group, since we
have equipped $\mathcal E$ with the discrete topology.  
To summarize we have shown:

\begin{theorem}
Let $(X,G)$ be a minimal dynamical system with abelian $G$ and $\mathcal E\subset \hat G$ the subgroup of (continuous) eigenvalues. Then the maximal equicontinuous factor is conjugate to $X \to \hat{\mathcal E}$, given by $x\mapsto j_x$, where $j_x:\mathcal E\to \T^1$ is defined by $j_x(\chi) = f_\chi(x)$, and the $G$-action on 
$\varphi\in\hat{\mathcal E}$ is given by
 $(t\cdot \varphi)(\chi) = \chi(t)\varphi(\chi)$.  

\end{theorem}
\newcommand{\EE}{H}
\begin{lemma}\label{lem-free} Let $\EE$ be a subgroup of $\hat G$.
Consider the $G$-action  $(t\cdot \varphi)(\chi) = \chi(t)\varphi(\chi)$ for $\varphi\in\hat \EE$ and $\chi\in \EE$.
This action is locally free if and only if $\hat G /\bar\EE$ is compact ($\bar\EE$ is the closure of $\EE$ in $\hat G$). 
In particular, the action is free if and only if $\EE$ is dense in $\hat G$.
\end{lemma}
\begin{proof}
Clearly $t\in G$ acts freely if and only if $t\cdot \varphi \neq \varphi$ for all $\varphi\in\hat \EE$, which is the case whenever $\chi(t) \neq 1$ for at least one $\chi\in\EE$. By continuity, the latter can be rephrased as $\chi(t) \neq 1$ for at least one $\chi\in\bar \EE$. 
Consider the exact sequence of abelian groups 
$$ 0 \to \widehat{\hat G/\bar \EE} \to G \stackrel{q}{\to} \hat{\bar\EE}\to 0$$
which is the dual to the exact sequence $0\to \bar\EE\to \hat G\to \hat G/\bar \EE\to 0$.

Suppose that $\hat G/\bar \EE$ is compact. This is the case whenever $ \widehat{\hat G/\bar \EE}$ is discrete.
Hence there exists $0\in U\subset G$, an open neighbourhood of the neutral element, such that $U\cap  \widehat{\hat G/\bar \EE} = \{0\}$. It follows that $q(U)\cong U$. Since for any $0\neq t\in \hat{\bar\EE}$, there exists
$\chi\in \EE$ such that $\chi(t)\neq 1$, we see that $U$ acts freely on $\hat \EE$. Now if $\EE$ is dense then $q$ is an isomorphism and we can take $U=G$. The converse, which we will not use, is left to the reader.
\end{proof}

We now consider the situation in which we have an additional homeomorphism $\Phi:X\to X$ which is compatible with the $G$-action in the sense that
$$ \Phi(t\cdot x) = \Lambda(t)\cdot \Phi(x)$$
for some group isomorphism $\Lambda:G\to G$. This will be relevant below when we consider substitution tilings.
If $f$ is an eigenfunction with eigenvalue $\chi\in\hat G$, then
$$ f(\Phi(t\cdot x)) = \chi(\Lambda(t)) f(\Phi(x)),$$
showing that $\Phi^*f$ is an eigenfunction with eigenvalue $\hat\Lambda \chi$. In particular, $\Phi^*$ preserves $\mathcal F$ and hence induces  an action $\Phimax$ on the maximal equicontinuous factor $\Xmax=\hat{\mathcal F}$ which is equivariant w.r.t.\ the $G$ action, and $\E$ is invariant under the dual isomorphism
$\hat\Lambda:\hat G\to \hat G$.
We ask the question: When is $\Phimax$  ergodic w.r.t.\ the Haar measure $\eta$?
This is precisely the case if the linear operator $U_\Phi$ defined on $L^2(\Xmax,\eta)$ by $U_\Phi\psi := 
\psi\circ \Phimax$ 
has, up to normalization, only one eigenfunction with eigenvalue $1$; that is, if $\psi = \psi\circ \Phimax$ implies that $\psi$ is constant. 

First note that if $f$ is a continuous eigenfunction with eigenvalue $\chi$, then $c_\chi := \frac{f(\Phi(x))}{f(x)}$ does not depend on $x$ and so the above equation reads
$$ \Phi^*f_\chi = c_{\hat\Lambda \chi} f_{\hat\Lambda \chi}$$
where $f_\chi$ is a normalized eigenfunction; i.e.,\ it satifies $f_\chi(x_0) = 1$ for some choice of $x_0$. 
Next note that normalized eigenfunctions form an orthogonal base in $L^2(\Xmax,\eta)$ and that 
$U_\Phi f_\chi = \Phi^*f_\chi$. We can therefore solve the eigenvalue equation $U_\Phi \psi = \psi$ in the above basis obtaining, for $\psi = \sum_{\chi\in\E} a_\chi f_\chi$, the equation
$$ a_{{\hat \Lambda}^{-1}\chi} =  c_{ \chi} a_{\chi}\quad \forall \chi\in \E.$$
Since $\sum_{\chi\in\E} |a_\chi|^2$ must be finite, and $|c_\chi|$  equals 1,
a solution different from $\psi=1$ can only exist if $\hat\Lambda$ admits a non-trivial periodic orbit in $\E$.
In particular we have:
\begin{lemma}\label{lem-ergodic}
If $G=\R^n$ and $\Lambda$ does not have any eigenvalues on the unit circle then $\Phi_{max}$ is ergodic with respect to $\eta$.
\end{lemma}
\begin{proof} In fact, upon identifying $\hat\R^n$ with ${\R^n}^*$,  $\hat\Lambda$ becomes the transpose and hence does not admit a non-trivial periodic orbit in $\hat G$ and therefore neither in $\E$.
\end{proof}
 \bigskip

We now consider conditions for pure point dynamical spectrum.
Let $\mu$ be a $G$-invariant Borel probability measure on $X$. 
An {\bf $L^2$-eigenfunction}
is an element $f\in  L^2(X,\mu)$ satifying the eigenvalue equation
(\ref{eq-eigen}). 
\begin{definition}
The measure dynamical system $(X,G,\mu)$ has pure point dynamical spectrum if
the $L^2$-eigenfunctions span $L^2(X,\mu)$. 
\end{definition}
Since $(X,G)$ is minimal, the maximal equicontinuous factor $(\Xmax,G)$ is also minimal, i.e.,\ $\Xmax$ is a group which contains $G/stab(X)$ as dense subgroup ($stab(X)$ being the stabilizer of $X$).
Hence a $G$-invariant Borel probability measure on $\Xmax$ is even $\Xmax$-invariant.\footnote{
Given that $\Xmax$ is metrisable, the $G$-invariant Borel probability measures are regular and therefore coincide via Riesz-representation theorem to $G$-invariant normed linear functionals on $C(X)$. By continuity of these functionals, $G$-invariance extends to $\Xmax$-invariance.}
Thus the only $G$-invariant Borel probability measure on $\Xmax$ is the Haar measure $\eta$.
\begin{theorem}\label{thm-1.3} Let $(X,G,\mu)$ be
a minimal dynamical system with ergodic Borel probability measure $\mu$. Assume that $(X,G,\mu)$ has finite minimal rank and that $\Ymax$ has full Haar measure. 
Then the following are equivalent:
\begin{enumerate}
\item
The continuous eigenfunctions generate
$L^2(X,\mu)$.
\item The minimal rank is $1$.
\item Proximality is a closed relation. 
\end{enumerate}
\end{theorem}
\begin{proof}
Equivalence of the last two assertions has already been shown.
The push forward of the measure $\mu$ under $\gmax$ is a $G$-invariant Borel probability measure and 
hence equals $\eta$.

Suppose that $mr=1$. By Lemma~\ref{lem-cr=mr} $cr=1$. Thus  the hypothesis that $\Ymax$ has full measure implies that $\gmax$ is a measure isomorphism and $\gmax^*(L^2(\Xmax,\eta))=L^2(X,\mu)$.
Since the linear span of the continuous eigenfunctions is dense in $L^2(\Xmax,\eta)$ and $\gmax$ is continuous, the linear span of the continuous eigenfunctions is 
dense in $L^2(X,\mu)$ as well.

Now let $mr=cr$ be greater than $1$ (but finite). Let $X':=\gmax^{-1}(\Ymax)$ and $\pi'$ be the restriction of $\gmax$ to $X'$.  Then $\pi':X'\to\Ymax$ is a closed map ($\gmax$ is closed) which is exactly $cr$-to-$1$ everywhere, and, for each $x\in\Ymax$, $\gmax^{-1}(x)=\{x_1,\ldots,x_{cr}\}$ with $d(x_i,x_j)\ge \delta_0$ for $i\ne j$ and $\delta_0>0$ as in Lemma~\ref{lem-14}. It follows that $\pi'$ must be injective on $\frac{\delta_0}2$ balls and thus $\pi':X'\to \Ymax$ a $cr$-to-$1$ local homeomorphism. It particular we can find an open set $U'\subset \Ymax$ such that ${\pi'}^{-1}(U')\subset \cup_{i=1}^{cr}B_{\delta_0/2}(x_i)$. $U'$ being open means of course that there exists an  open set $U\subset \Xmax$ such that $U'=U\cap\Ymax$. 
Let $U_i:=B_{\delta_0/2}(x_i)\cap {\gmax}^{-1}(U)$. Since $\Ymax$ has full measure also $X'$ has full measure and thus $\gmax\left|_{U_i}\right.:U_i\to U$ is a function which is almost everywhere bijective and bi-measurable.

For each $U,U_i$, as above, let $\eta_U^i$ be the push forward of $\eta|_U$ by $(\gmax|_{U_i})^{-1}$. That is, $\eta_U^i(A):=\eta(\gmax(A))$ for Borel $A\subset U_i$. Then $\mu|_{U_i}<<\eta_U^i$: let $J_{U_i}:U_i\to\R$ be the Radon-Nikodym derivative of $\mu|_{U_i}$ with respect to $\eta_U^i$. By uniqueness of Radon-Nikodym derivative, the $J_{U_i}$ paste together to give
a Borel-measurable function
$J:X\to\R$. It follows from the $G$-invariance of $\mu$ and $\eta$ that $J$ is also $G$-invariant. By ergodicity of $\mu$, $J$ must be $\mu$-almost everywhere constant and this constant therefore equal to $1/cr$.

Now let $h:=\mathbf{1}_{U_1}-\mathbf{1}_{U_2}$ be the difference of indicator functions for some $U,U_1,U_2$ as above, and suppose that $f:X\to\mathbb{C}$ is a continuous eigenfunction. Then there is $f':\Xmax\to\mathbb{C}$ so that 
$f|_{U_i}=f'\circ\gmax|_{U_i}$ for $i=1,2$. The scalar product between $f$ and $h$ is
$ \langle f,h\rangle=\int_{U_1}f\,d\mu-\int_{U_2}f\,d\mu=\int_Uf'\,\frac{1}{cr}d\eta-\int_Uf'\,\frac{1}{cr}d\eta=0$. That is, $h$ is orthogonal to all continuous eigenfunctions. As $\mu(U_1)=
\frac{1}{cr}\eta(U)$ and $\eta(U)>0$ (since $U$ is open),
$h$ is not the zero function, and we see that the linear span of the continuous eigenfunctions is not dense in $L^2(X,\mu)$.

\end{proof}
\begin{cor}
Consider a minimal dynamical system with ergodic $G$-invariant Borel probability measure $(X,G,\mu)$ 
with finite minimal rank whose dynamical eigenfunctions are all
continuous. Assume that $\Ymax$ has full 
Haar measure. 
Then the dynamical spectrum 
of $(X,G,\mu)$ is pure point if and only if the proximality
relation is closed. 
\end{cor}

\renewcommand{\L}{\mathcal L}
\section{Proximality for tilings and Delone sets}
\subsection{Preliminaries}
By a {\bf tile} $\tau$ in $\R^n$ we mean here a compact subset of $\R^n$ which is the closure of its interior.
It is sometimes useful to decorate tiles with marks and then one should rather speak of a tile as an ordered pair $\tau=(spt(\tau),m)$ where $spt(\tau)$, the {\bf support} of $\tau$ is compact and the closure of its interior
and $m$  a {\bf mark} taken from some finite set of marks. The {\bf interior} of a tile is then simpliy the interior of its support: $\mathring{\tau}:=int(spt(\tau))$.
Two tiles $\tau=(spt(\tau),m)$ and $\sigma$ are {\bf translationally equivalent} if there is a $v\in\R^n$ with $\tau+v:=((spt(\tau)+v,m)=\sigma$. 

A {\bf patch} is a collection of tiles with pairwise disjoint interiors, the {\bf support} of a patch $P$, $spt(P)$, is the union of the supports of its constituent tiles, the {\bf diameter} of $P$, $diam(P)$, is the diameter of its support, and a {\bf tiling} of $\R^n$ is a patch with support $\R^n$. We denote the translation action on patches (and tilings) also by $P\mapsto P-v$, $v\in\mathbb R^n$.  

A collection $\Omega$ of tilings of $\R^n$
has {\bf translationally finite local complexity} (FLC) if it is the case that for each $R$ there are 
only finitely many translational equivalence classes of patches $P\subset T\in\Omega$ with $diam(P)\le R$. 
It is very useful to consider a metric topology on sets of tilings. This is expressed with the help of $R$-patches. 
Given a tiling $T$ and $R\geq 0$ the patch $B_R[T]:=\{\tau\in T:\bar B_R(0)\cap spt(\tau)\ne\emptyset\}$ is called
the {\bf $R$-patch} of $T$ at $0$. 
If $\Omega$ is a collection of tilings of $\R^n$ with FLC, then the following metric $d$ can be used:
 \begin{equation}\label{metric}
d (T, T') := \inf \{\frac\epsilon{\epsilon+1} : \exists \|v\|,\|v'\|\leq \frac\epsilon2:
B_{\frac1\epsilon} [T-v]=B_{\frac1\epsilon}[T'-v']\}.
\end{equation}
In other words, in this metric two tilings are close if a small translate of one agrees with the other in a large neighborhood of the origin.  If we take out the possibility of translation by a small vector 
we get another metric:
\begin{equation}\label{metric0}
d_0 (T, T') := \inf \{\frac\epsilon{\epsilon+1} : 
B_{\frac1\epsilon} [T]=B_{\frac1\epsilon}[T']\}.
\end{equation}
which does not induce the same topology  but is also useful. 

We will call a collection $\Omega$ of tilings of $\R^n$ an {\bf n-dimensional tiling space} if $\Omega$ has FLC, is closed under translation ($T\in\Omega$ and $v\in\R^n\Rightarrow T-v\in\Omega$), and is compact in the metric $d$. (All tiling spaces in this article are assumed to have FLC, but we will occasionally include the FLC hypothesis for emphasis.) For example, if $T$ is an FLC tiling of $\R^n$, then 
$$\Omega_T=\{T': T' \mbox{ is a tiling of $\R^n$ and every patch of $T'$ is a translate of a patch of $T$}\}$$ 
is an $n$-dimensional tiling space, called the {\bf hull} of $T$ (\cite{AP}).

A $n$-dimensional tiling space $\Omega$ is {\bf repetitive} if for each patch $P$ with compact support that occurs in some tiling in $\Omega$ there is an $R$ so that for all $T\in\Omega$ 
a translate of $P$ occurs as a sub-patch in $B_R[T]$. 
If $\Omega$ is repetitive, then the action of $\R^n$ on $\Omega$ by translation is minimal.

For a tiling $T$ let $p:T\to\R^n$ satisfy: $p(\tau)\in spt(\tau)$ and, if $\tau,\tau+v\in T$, $p(\tau+v)=p(\tau)+v$. Such an assignment $p$ will be called a {\bf puncture map}. If the tiling has FLC then
the set of its punctures $p(T)=\{p(\tau):\tau\in T\}$ is a Delone set, i.e.,\ of subset of $\R^n$ which is  uniformly discrete and  relatively dense\footnote{Uniformly discrete means there is $r>0$ so that $\sharp B_r(x)\cap\L\le1$ for all $x\in\R^n$ and relatively dense means that 
there is $R$ so that $B_R(v)\cap\L\ne\emptyset$ for each $x\in\R^n$.}.
Notice that, if $p$ and $p'$ are two choices of puncture maps for an FLC tiling
then there is a finite set $F$ such that $p(T)-p(T)\subset p'(T)-p'(T)+F$. 

A puncture map defines a transversal in the hull $\Omega_T$ of a tiling $T$, namely the set
of $T'\in\Omega_T$ such that $0\in p(T)$. Restricting the metrics $d$ and $d_0$ defined above to this 
transversal they become equivalent. 

The definitions we have made for tilings all have analogs for Delone sets and whether we deal with tilings or Delone sets is mainly a matter of convenience. One could, for instance, 
represent a tiling $T$ by the Delone set of its punctures\footnote{Strcitly speaking, one might have to consider the marked Delone set (or Delone multi-set) $\{(p(\tau),m):\tau\in T, m=m(\tau),\text{ the mark of }\tau\}$ for that. Everything we do below with Delone sets could as well be done with marked Delone sets - we trust the interested reader to make the necessary adjustments.}. But there are other possibilities. If the tiling is polyhedral (and FLC), one could represent it also by the Delone set of its vertices.
 On the other hand, when dealing with Delone sets we can carry over the notions defined
above for tilings if we consider the Dirichlet tiling associated with the Delone set (the tiling defined by the dual of the Voronoi complex) which is a polyhedral tiling having the points of the Delone set as vertices.
With this in mind, we may define $R$-patches for a Delon sets as the $R$-patches of its associated Dirichlet tiling
and hence get a metric on collections of Delone sets as well. Alternatively one could  use the more standard definition of the $R$-patch at $x$ of the Delone set $\L$ as the set $\{y\in\L:\|x-y\|\leq R\}$ to obtain a metric.
The two metrics are different but define the same topology and all our results are independent of  choice of metric. 
 
\subsection{Strong proximality of tilings and the Meyer property}

We consider dynamical systems $(\Omega,\R^n)$, where $\Omega$ is an $n$-dimensional tiling space and 
$\R^n$ acts by translation. The following two definitions and results are mostly formulated for tilings but have obvious counterparts for Delone sets. 

\begin{definition}[strong proximality] Two tilings $T_1,T_2\in\Omega$ are called strongly proximal if, for
all $r$, there exists $v\in \R^n$ 
such that $B_r[T_1-v] = B_r[T_2-v]$.
Strong proximality is thus proximality for the metric $d_0$. 
\end{definition}
\begin{definition}[strong regional proximality] Two tilings $T_1,T_2\in\Omega$ are called strongly regionally proximal if, for
all $r$, there exist $S_1,S_2\in \Omega_T$ and $v\in \R^n$ such that 
$$ B_r[T_1] = B_r[S_1],\quad B_r[T_2]=B_r[S_2],\quad B_r[S_1-v] = B_r[S_2-v].$$
Strong regional proximality is thus regional proximality for the metric $d_0$. 
\end{definition}
An important question is: For what
classes of tilings does proximality (regional proximality) imply strong proximality (strong regional proximality)? If the relations are the same then proximality (regional proximality) becomes a purely combinatorial
property. 

Recall that a {\bf Meyer set} is a Delone set $\L$ such that $\L-\L$ is uniformly discrete. In particular, a Meyer set is always FLC. We will say that a tiling $T$ is a {\bf Meyer tiling}, or {\bf has the Meyer property}, if it has FLC and the image $\L=p(T)$ of a puncture map is a Meyer set. Although the set $\L$ depends on the puncture map, the property of its being Meyer does not.
If $\L$ is Meyer then so also is $\L-\L$: 
\begin{prop}[Meyer]
For a Meyer set $\L$ one has that all finite combinations $\L \pm \L \pm\cdots \pm\L$ (any choice of signs) are also Meyer.
\end{prop} 
A proof and various different characterizations of Meyer sets can be found in \cite{M}. 

Note that, if $\L'$ is any element in the hull $\Omega_\L$ of a Delone set $\L$ then $\L'-v'\in\L-\L$ for all $v'\in\L'$.
Indeed, for all $r$ exists $v$ such that $B_r[\L']=B_r[\L-v]$. Take $v'\in\L'$. Then $0\in\L-v-v'$.
Hence $v\in \L-v'$ and so each point of $\L'$ lies in $\L-\L+v'$. It follows that  $\L_1 \pm \L_2 \pm\cdots \pm\L_k$ (for any choice of signs) is Meyer for all $\L_i\in\Omega_\L$.

Now suppose that $\L_1,\L_2\in\Omega_\L$ satisfy $\L_1\cap\L_2\neq \emptyset$. Take $v\in \L_1\cap\L_2$.
Then $\L_1-\L_2=(\L_1-v)-(\L_2-v)\subset (\L-\L)-(\L-\L)$. Hence $\bigcup_{\L_1,\L_2\in\Omega_{\L}:\L_1\cap\L_2\neq\emptyset}\L_1-\L_2$ is contained in a uniformly discrete set and thus cannot contain an accumulation point.
 
\begin{theorem} \label{FLCproximality}
Let $\Omega$ be the hull of a repetitive 
Meyer tiling. The proximal relation (regional proximal relation) on $\Omega$ coincides
with the strong proximal relation (strong regional proximal relation).
\end{theorem}
\begin{proof}
We provide the proof for the regional proximality relation. The proof for proximality is similar and a bit simpler. Let $(T,T')\in \Qr$. Hence for all $r>0$ there exist $t,t',v,z\in \R^n$ with $|z|\leq \frac1r$ such that:
\begin{enumerate}
\item $B_r[T] = B_r[T-t]$,
\item $B_r[T'-z] = B_r[T-t']$, and
\item $B_r[T-t-v] = B_r[T-t'-v]$.
\end{enumerate}
In fact, since $P_\epsilon$ is open and the orbit of $T$ dense we can take the $x'$ and $y'$ required in the definition of regional proximality to lie in that orbit and thus be of the form $T-t$ and $T-t'$. 
Our aim is to show that we can take $z=0$. Let $\L=p(T)$ and $\L'=p(T')$ where $p$ is a puncture map such that $0\in\L$.
(1) implies that $t\in \L$. (3) implies that $(t'+v)-(t+v)\in \L-\L$. Hence $t' \in \L-\L+\L$. Let $a\in \L-t'$, $|a|<r$.
Hence $a\in \L - \L + \L -\L$. 
Now (2) implies that $a\in \L'-z$. Hence $z\in 
\L' - \L + \L - \L +\L$ showing that $z$ takes values in a uniformly discrete set. Thus, since $|z|\leq \frac1r$, there exists $r_0$ such that $z=0$ if $r\geq r_0$.     
\end{proof}
\begin{corollary}\label{finitely many pairs} Consider the regional proximal relation $\Qr$ on the hull
of a repetitive Meyer tiling and let $s>0$. Up to translation, there are only finitely many
pairs of patches of the form $(B_s[T],B_s[T'])$ with $(T,T')\in \Qr$.
\end{corollary}
\begin{proof}
We denote by $[(P,P')]$ the translational congruence class of a pair of patches $(P,P')$. Let $(T,T')\in\Qr$ and $s>0$. By FLC there is a finite list $\{P_1,\cdots,P_k\}$ of $s$-patches such that any $s$-patch is a translate of some $P_i$. Hence the set of translational congruence classes 
$[(B_s[S],B_s[S'])]$ with $(T,T')\in\Qr$ is of the form
$\{[(P_i,P_j-t_{i,j})]:(i,j,t_{i,j})\in I\}$ with some subset $I$ of $\{1,\cdots,k\}^2\times B_{s'}(0)$ and some finite $s'$. 
We need to show that $I$ is finite. Assume the contrary. Let $(i,j,t)$ be an accumulation point of $I$ and let $v_{i,j}\in p(P_i)-p(P_j)$, $p$ a puncture map.
 Then $t+v_{i,j}$ is an accumulation point of $\bigcup_{(i,j,t_{i,j})\in I}p( P_i)-(p(P_j)-t_{i,j})$.

By the last theorem there exists $v\in \R^n$ and $S,S'\in\Omega$ such that $B_s[T]=B_s[S]$,
 $B_s[T']=B_s[S']$ and $B_s[S-v]=B_s[S'-v]$. 
 It follows that the set $\{[(B_s[T],B_s[T'])] : (T,T')\in \Qr\}$ is contained in the set
 $\{[(B_s[S],B_s[S'])] : p(S)\cap  p(S')\neq\emptyset\}$. The latter set thus contains
 $\{[(P_i,P_j-t_{i,j})] : (i,j,t_{i,j})\in I\}$ and hence $t+v_{i,j}$ must be
 an accumulation point of $\bigcup_{S,S'\in\Omega:p(S)\cap p(S')\neq\emptyset} p(S)-p(S')$.  
 But, as we saw in the discussion preceding Theorem \ref{FLCproximality}, that set does not contain any accumulation points.
 \end{proof}
If $S=\{T_i\}_i$ is a set of tilings or Delone sets,  we denote by $B_R[S]$ the corresponding collection of $R$-patches, $B_R[S] := \{B_R[T_i]\}_i$.
\begin{corollary}\label{cor-bound} Consider the dynamical system of a
repetitive 
Meyer tiling. There is an $R_0$ such that for all $y\in\Xmax$ and $R\geq R_0$
$$ \sup\{l: \#B_R[\gmax^{-1}(y)-v]\geq l\text{ for all } v\in\R^n\}=cr .$$
\end{corollary}
\begin{proof}
We suppose that the minimal rank is finite.
By Lemma~\ref{lem-14} we have that for all $y$: 
\begin{eqnarray*}
cr & = & \sup\{l: \exists T_1,\ldots,T_l\in\gmax^{-1}(y) \text{ such that }\forall i\neq j \text{ and } \forall v,\,d(T_i-v,T_j-v)\geq\delta_0\} \\
&=& \sup\{l: \exists T_1,\ldots,T_l\in\gmax^{-1}(y) \text{ s.\ th. }\forall i\neq j,\, \forall R\geq R_0,\forall v,\,B_R[T_i-v]\neq B_R[T_j-v]\}, \end{eqnarray*}
where we have used the Meyer property for the last equation and taken $R_0$ suitably large depending on $\delta_0$. 
The case in which the minimal rank is not finite is easily obtained.
\end{proof}

Let $R>0$ and
$$n^R(x):=\#B_R[\gmax^{-1}(x)] = \#\{B_R[T]:T\in \gmax^{-1}(x)\}.$$
For later use we show
\begin{lemma} \label{lem-upper} For Meyer tilings with finite maximal rank
$n^R$ is upper semi-continuous, i.e.,\ $\{x:n^R(x)\ge k\}$ is closed for all $k$. In particular,  if $R$ is sufficiently large (depending on the $\delta_0$ of Lemma \ref{lem-14}) then
$D^R:=\{x\in \Xmax:n^R(x) = cr\}$ has strictly positive measure.
\end{lemma}
\begin{proof}
Recall that, for $T\in\Omega$, $B_R[T]$
is defined to be the collection of tiles in $T$ that meet the closed
(rather than open) ball of radius $R$ at $0$. It hence follows that, for all $R$, there exists $R'>R$ such that
$B_{R'}[T]=B_{R}[T]$.  

The fibers of $\pi$ satisfy the following property: If $(x_n)_n$ is a sequence in $\Xmax$ tending to $x$ then
all accumulation points of sequences $(T_n)_n$, $T_n\in \gmax^{-1}(x_n)$ are contained in $\gmax^{-1}(x)$.
Since $\gmax^{-1}$ is finite it is uniformly discrete and hence if $n$ is large enough we can define a map $f_n:\gmax^{-1}(x_n)\to \gmax^{-1}(x)$ by saying that $f_n(T)$ is the element of $\gmax^{-1}(x)$ which is closest to $T$.

Now let $B_R[f_n(T)] = B_R[f_n(T')]$. By the first remark there exists $R'>R$ such that $B_{R'}[f_n(T)] = B_{R'}[f_n(T')]$. We may assume that $n$ is large enough so that $d(f_n(T),T)$ 
is small enough to guarantee that
$B_{R'}[T-v] = B_{R'}[f_n(T)]$ for some $|v|\leq R'-R$. Likewise, we may assume that $n$ is large enough 
so that $B_{R'}[T'-v'] = B_{R'}[f_n(T')]$ for some $|v'|\leq R'-R$. We may suppose that $R'-R$ is small so that $B_{R'}[f_n(T)] = B_{R'}[f_n(T')]$ and
 the Meyer property imply that $v=v'$ and therefore $B_{R'}[T-v] =B_{R'}[T'-v] $. The latter implies
$B_{R}[T] =B_{R}[T'] $. This shows that $n^R(x)\geq n^R(x_n)$ and hence upper semi-continuity of $n^R$.

We have $\#\gmax^{-1}(x) \geq n^{R'}(x) \geq n^{R}(x)\geq cr$ if $R'>R(\delta_0)$
where the $\delta_0$ is the one from  Lemma~\ref{lem-14}. Thus
$$D^R:=\{x\in \Xmax:n^R(x) = cr\}=\{x\in \Xmax:n^R(x)\leq cr\}$$ showing that $D^R$ is
is open. Since $D^R\ne\emptyset$, the Haar measure of $D^R$ is positive.
\end{proof}

\section{Proximality in model sets}
\newcommand{\Torus}{\mathbb T}
Cut \& project patterns, or model sets, are characterized by the way they are constructed. We outline the construction here refering the reader to \cite{FHK,M,BLM} for a thorough desription.
 
The defining data of a model set are a lattice (co-compact subgroup) $\Gamma\subset \R^n\times H$ in the product of $\R^n$ with a locally compact abelian group $H$ such that $\R^n$ is in irrational position w.r.t.\ $\Gamma$, and a window $K$ (or acceptance domain, or atomic surface) which is a compact subset of $H$. 
We denote the boundary of $K$ by $\partial K$ and the quotient group $\R^n\times H/\Gamma$ by $\mathbb{T}$.
The latter is a compact abelian group which is often referred to as the LI-torus.
$\R^n$ acts on $\Torus$ by rotation: $v\cdot ((w,h) + \Gamma)=(w+v,h)+ \Gamma$. 
$(\Torus,\R^n)$ is thus an equicontinuous dynamical system.

Let $\pi:\R^n\times H\to \R^n$ be the projection onto the first factor and
$\pi^\perp:\R^n\times H\to H$ be the projection onto the second factor.
We make the following standard assumptions.  
\begin{itemize}
\item  The restrictions of $\pi^\|$ and $\pi^\perp$ to
$\Gamma$ are one to one,
\item  the restrictions of $\pi^\|$ and $\pi^\perp$ to
$\Gamma$ have dense image,
\item $K$ is the closure of its interior,
\item the stabiliser of $K$ in $H$ is trivial; that is, $h+K=K$ implies $h=0$. 
\end{itemize}
The data $(\R^n,H,\Gamma,K)$ 
determine a whole family of point patterns in $\R^n$ which we identify with $\R^n\times\{e\}$, $e\in H$ the neutral element. Indeed, for $x\in \R^n \times H$ let
$$\MS_x:=\R^n\times \{e\} \cap(\Gamma+x-\{0\}\times K).$$
It is well known that under the above assumptions, $\MS_x$ is a repetitive Delone set.
A Delone set arising in this way is called a {\bf model set}. Furthermore,
$\MS_x=\MS_y$ if and only if 
$x-y\in\Gamma$. We define the set $S$ of {\bf singular points} by 
$$ S :=\{ x\in \R^n\times H : \pi^\perp(x)\in \partial K+\pi^\perp(\Gamma)\} = \R^n\times\{e\} + \Gamma + \{0\}\times\partial K$$ and denote by $N\!S$ its complement. 
\begin{prop}
$N\!S$ is a dense $G_\delta$ set. In particular it is non-empty.
\end{prop}
\begin{proof} By our assumptions $\partial K$, and therefore also $R^n \times\{e\}  + \{0\}\times\partial K$, has empty interior and hence so does $S$. This is 
a simple application of the Baire Category Theorem; see \cite{Schlottmann} or \cite{FHK} for the case $H=\R^k$.
\end{proof}
Let us suppose for simplicity that $0\notin S$ and consider
the hull  $\Omega_\MS$ of $\MS=\MS_0$. 
It is well-known that $\MS_x$ and $\MS_y$ are locally indistiguishable provided $x,y\in N\!S$. 
Furthermore $\MS_x=\MS_y$ if and only if $x-y\in\Gamma$. 
Thus $\Omega_\MS$
is the completion of the set $NS/\Gamma$ w.r.t.\ the metric $\delta(x+ \Gamma,y+\Gamma) = d(\MS_x,\MS_y)$.The metric $\delta$ does not extend continuously in the (quotient of the product) topology of $\Torus$ but the converse is  
the basis of one of the main structural theorems for model sets.
\begin{theorem}
The map $\{\MS_y\in\Omega_\MS:y\in NS\} \ni \MS_x\mapsto x +\Gamma \in\Torus$ extends to a continuous surjection
$$\mu:\Omega_\MS\to\Torus$$ which is equivariant with respect to the $\R^n$-action and one-to-one precisely on $N\!S/\Gamma$, i.e.,\ precisely the non-singular points have a unique pre-image. 
\end{theorem}
\begin{proof}
In the context in which $H$ is a real vector space a proof can be found in \cite{FHK}. 
The case of more general groups $H$ can be found in \cite{BLM} .
\end{proof}
\begin{cor} Model sets have minimal rank $1$. In particular, 
$(\Torus,\R^n)$ is the maximal equicontinuous factor, $\mu=\gmax$, and
two elements are proximal if and only if they are mapped to the same point by $\mu$.  
\end{cor}
\begin{proof} The set of fiber distal points includes $N\!S/\Gamma$ which is
  non empty. Hence $cr=mr=1$. By Lemma~\ref{lem-cr=mr} the proximality relation
  coincides with the equicontinuous structure relation. 
\end{proof}
A model set is called {\bf regular} if $\partial K$ has $0$ measure (w.r.t.\  Haar measure on $H$).  
\begin{theorem}
For a regular model set, $\Ymax$ has full Haar-measure. Moreover,
if the model set is not regular, then $\Ymax$ has $0$ Haar measure.
\end{theorem}
\begin{proof} 
Since the proximality relation  coincides with the equicontinuous structure relation,
the set of distal points coincides with the non-singular points and 
$\Ymax=N\!S/\Gamma$. Clearly $\partial K$ has
strictly positive measure if and only if the complement of  $\Ymax=N\!S/\Gamma$
has strictly positive Haar measure. 
By ergodicity 
of the Haar measure $S/\Gamma$ must have full Haar measure if its
measure is strictly positive.
\end{proof}

\section{Proximality for Meyer substitution tilings}\label{preliminaries}

\subsection{Basic notions}

Suppose that $\mathcal{A}=\{\rho_1,\ldots,\rho_k\}$ is a set of translationally inequivalent tiles (called {\bf prototiles}) in $\R^n$ and $\Lambda$ is an expanding linear isomorphism of $\R^n$, that is, all eigenvalues of $\Lambda$ have modulus strictly greater than $1$. 
A {\bf substitution} on $\mathcal{A}$ with expansion $\Lambda$ is a function $\Phi:\mathcal{A}\to\{P:P$ is a patch in $\R^n\}$
with the properties that, for each $i\in\{1,\ldots,k\}$, every tile in $\Phi(\rho_i)$ is a translate of an element of $\mathcal{A}$, and $spt(\Phi(\rho_i))=\Lambda(spt(\rho_i))$. Such a substitution naturally extends to patches whose elements are translates of prototiles by $\Phi(\{\rho_{i(j)}+v_j:j\in J\}):=\cup_{j\in J}(\Phi(\rho_{i(j)})+\Lambda v_j)$. A patch $P$ is {\bf allowed} for $\Phi$ if there is an $m\ge1$, an $i\in\{1,\ldots,k\}$, and a $v\in\R^n$, with $P\subset \Phi^m(\rho_i)-v$. The {\bf substitution tiling space} associated with 
$\Phi$ is the collection $\OP:=\{T:T$ is a tiling of $\R^n$ and every finite patch in $T$ is allowed for $\Phi\}$. Clearly, translation preserves allowed patches, so $\R^n$ acts on $\OP$ by translation.

The substitution $\Phi$ is {\bf primitive} if for each pair $\{\rho_i,\rho_j\}$ of prototiles there is a $k\in\N$ so that a translate of $\rho_i$ occurs in $\Phi^k(\rho_j)$.  If $\Phi$ is primitive then $\OP$ is repetitive.

If the translation action on $\Omega$ is free (i.e., $T-v=T\Rightarrow v=0$), $\Omega$ is said to be {\bf non-periodic}. If $\Phi$ is primitive and $\OP$ is FLC and non-periodic then $\OP$ is compact in the metric described above, $\Phi:\OP\to\OP$ is a homeomorphism, and the translation action on $\OP$ is minimal and uniquely ergodic ( \cite{AP}, \cite{S3}, \cite{sol}). In particular, $\OP=\Omega_T$ for any $T\in\OP$. It will be with respect to the unique ergodic measure $\mu$ on $\OP$ when we speak about the dynamical spectrum and $L^2$-eigenfunctions. Note that $\Phi$ preserves regional proximality so there is an induced homeomorphism $\Phimax$ on the maximal equicontinuous factor $\Xmax$ of $\OP$.

All substitutions will be assumed to be primitive, aperiodic and FLC.

\begin{theorem}[\cite{S}] \label{thm-ev1}
All $L^2$-eigenfunctions of a substitution tiling space can be chosen to be continuous.  
\end{theorem}

We call a substitution a {\bf Meyer substitution} if every tiling $T\in \OP$ has the Meyer property
(that is, the set of punctures $p(T)$ is a Meyer set).
This does not depend on the choice of punctures and hence holds true also if punctures are control points in the sense of \cite{LS}. We consider below primitive aperiodic Meyer substitutions. In this context, $\OP$ is minimal and the Meyer property is satisfied for all $T\in \OP$ if it is satisfied for a single one.
\begin{definition}[Meyer substitution tiling]
A Meyer substitution tiling is a tiling in the hull of a
primitive aperiodic Meyer substitution.
\end{definition}

\subsection{Finite rank and fiber distality of Meyer substitutions}
Recall that the maximal rank of a tiling is
$$ \sup\{\#\gmax^{-1}(x):x\in\Xmax\} $$
which, of course, bounds the minimal rank.

\begin{theorem}\label{finite-to-one} A Meyer substitution tiling system
has finite maximal rank.  
\end{theorem}

\begin{proof}
By 
Corollary~\ref{finitely many pairs} there is $N$ so that
$\#\{B_1[T']:\gmax(T')=\gmax(T)\}\le N$ for all $T\in \Omega$. Suppose
that $T_1,\ldots,T_m$ are distinct tilings in $\Omega$ with 
$\gmax(T_i)=\gmax(T_j)$ for all $i,j$. Let $R$ be large enough so that
$B_R[T_i]\ne B_R[T_j]$ for all $i\ne j$ and let $k$ be large enough so
that $\Lambda^k(B_1(0))\supset B_R(0)$. Then $B_1[\Phi^{-k}(T_i)]\ne 
B_1[\Phi^{-k}(T_j)]$ for $i\ne j$. Thus $m\le N$ and $\gmax$ is at
most $N$-to-$1$. 
\end{proof}

The following extends the one-dimensional result of \cite{BK}.

\begin{theorem}\label{cr} For a Meyer substitution tiling system
$\Ymax$ has full Haar measure.
\end{theorem}
\begin{proof} Let $\Phi$ be a Meyer subtitution with tiling space $\Omega$ and with linear expansion  $\Lambda$.
There is $k>0$ with $\Lambda^k B_R(0)\supset B_R(0)$ and hence $B_R[\Phi^k(T)]=B_R[\Phi^k(B_R[T])]$ for all $T\in\Omega$. It follows from this, and $\Phimax\circ \gmax=\gmax\circ \Phi$, that $n^R(\Phimax^k(x))\le n^R(x)$.  Hence $D^R=\{x\in \Xmax:n^R(x) = cr\}$ is 
invariant under $\Phimax^k$. 
By Lemma~\ref{lem-upper}, $D^R$ has strictly positive measure. By Lemma~\ref{lem-ergodic}, 
 $\Phimax$ is ergodic with respect to Haar measure and
thus $\eta(D^R)=1$. Using Lemma \ref{lem-cr=mr},  $
\eta(\Ymax)=\eta(\{x:\sharp \gmax^{-1}(x)=cr\})=\eta(\cap_{R>0}D^R)=1$.
\end{proof}

\begin{corollary} The tiling flow on a Meyer substitution tiling space 
  has pure discrete spectrum if and only if the
  proximal relation is closed. 
\end{corollary}

\subsection{Pisot family substitutions}
Geometry places rather strong conditions on the collection, spec($\Lambda$), of eigenvalues of the expansion matrix $\Lambda$ of a substitution. To begin with, all elements of spec($\Lambda$) must have absolute value greater than 1, simply because $\Lambda$ is an expansion. To say more, it is convenient to introduce the notion of a family. Let $p$ be a monic integer polynomial and let $c>0$ be a real number. A collection of complex numbers of the form $$F_{p,c}:=\{\lambda\in\mathbb C:p(\lambda)=0, |\lambda|\ge c\}$$  is called a {\bf family}. That is, a family is a collection of all the algebraic conjugates of some algebraic integer $\lambda$ that have absolute value at least as great as that of $\lambda$. In general, the elements of spec($\Lambda$) must be algebraic integers (\cite{K2},\cite{LS2}), and if $\Lambda$ is diagonalizable over $\mathbb{C}$, spec($\Lambda$) must be a union of families (\cite{KS}). In the special case that $\Lambda=\lambda I$ the substitution is called {\bf self-similar} and the tiling space has a nontrivial equicontinuous factor
(equivalently, the $\R^n$-action has eigenvalues) if and only if $\lambda$ is a Pisot number - an algebraic integer greater than 1, all of whose algebraic conjugates have absolute value less than 1. Let us call a family a $F_{p,c}$ a {\bf Pisot family} if $c=1$ and no element of $F_{p,c}$ has absolute value 1.

Recall that eigenvalues for a group action are continuous characters. When the group is $\R^n$,
any such character takes the form $\chi(x)=e^{2\pi i\langle x,\beta\rangle}$ for some $\beta\in\R^n$. In this context, it is customary, as in the following theorem, to call $\beta$ (rather than $\chi$) an eigenvalue of the action.
\begin{theorem}[\cite{LS}]\label{thm-Pisot-Meyer}
Consider a primitive FLC substitution with diagonalisable expansion matrix $\Lambda$. Suppose that $spec (\Lambda)$ consists of algebraic conjugates with the same muliplicity. 
The following are equivalent:
\begin{enumerate}
\item the substitution is Meyer, 
\item spec($\Lambda$) is a Pisot family,
\item the eigenvalues are relatively dense in $\mathbb R^n$, 
\item the maximal equicontinuous factor is non-trivial.
\end{enumerate}
\end{theorem}

To capture all these desirable qualities, we say that a substitution is a {\bf Pisot family substitution} if it is primitive, aperiodic, FLC, and its linear expansion is diagonalizable over $\mathbb{C}$ and has a Pisot family spectrum with all elements of the same multiplicity. The {\bf degree} of such a substitution is the algebraic degree of the elements of the Pisot family, and its {\bf multiplicity} is the common multiplicity of those elements as eigenvalues of the expansion. 
\bigskip

The above theorem shows that, under the stated assumptions, the maximal equicontinuous factor is either a single point, or the set of eigenvalues contains a sub-group which is relatively dense in $\hat\R^n$. We will improve this result by determining completely the form of the group of eigenvalues and the maximal equicontinuous factor for Pisot family substitutions. The key result is an extension of Cor.~\ref{finitely many pairs}.

\begin{theorem}\label{thm-thm1} Let $\Omega$ be the continuous hull of an aperiodic primitive FLC substitution $\Phi$ with expansion matrix $\Lambda$. Let $g:\Omega\to \mathbb{X}$ be a factor of the tiling flow such that the
$\mathbb{R}^n$-action on $\mathbb{X}$ is locally free and
\begin{itemize}
\item[(1)] $g(T)=g(T')\implies g(\Phi(T))=g(\Phi(T')).$
\end{itemize} 
Let $p$ be a puncture map.
Then the set $$\mathcal{O}:=\{x: \exists T,T' \mbox{ with }g(T)=g(T'),
\tau\in T, \tau'\in T' ,\mathring{\tau}\cap\mathring{\tau}'\ne \emptyset, \mbox{ such that } p(\tau)-p(\tau')=x\}$$ 
is finite. 
\end{theorem}

\begin{proof}
Suppose $\mathcal{O}$ is infinite. Then there are sequences $(T_n)_n,(T_n')_n\in \Omega$ with $\tau\in T_n, \tau_n'\in T_n'$ with $g(T_n)=g(T_n')$ for all $n$,
$\mathring{\tau}\cap\mathring{\tau_n'}\ne\emptyset$ and $p(\tau_m')\ne p(\tau_n')$ for $m\ne n$. 
The conclusion of the theorem is independent of the choice of puncture, so we may assume that the puncture of any tile is in the interior of the support of the tile.
Translating, we may also assume $0=p(\tau)$.
Passing to a subsequence we may assume that:
\begin{itemize}
\item[(2)] { $T_n\to T, T_n'\to T'$} for some $T,T'$;
\item[(3)] {$p(\tau_n')\to p(\tau'),\tau'\in T'$}; and
\item[(4)] $\Phi^n(T)\to \bar{T}$ and $\Phi^n(T')\to
\bar{T}'$ for some $\bar T,\bar T'$.
\end{itemize}
Letting $x_n:=p(\tau'_n)-p(\tau')$, there is a function $m:\N\to\N$, with $m(n)\to\infty$, and an $x\ne0$, but close enough to 0 so that it acts freely, such that
\begin{itemize}
\item[(5)] $\Lambda^{m(n)}(x_n)\to x$. 
\end{itemize}
Note that since $0\in \mathring{\tau}$ and $\tau\in T_n$, 
\begin{itemize}
\item[(6)] $\Phi^{m(n)}(T_n)\to \bar{T}$
\end{itemize}
and likewise, since, at least for large $n$, $T_n'-x_n$ and $T'$ have exactly the same tiles at the origin, 
\begin{itemize}
\item[(7)] $\Phi^{m(n)}(T_n'-x_n)\to \bar{T}'$. 
\end{itemize}
By continuity of $g$, $g(T_n)=g(T_n')$ implies $g(T)=g(T')$, and hence (1),(4) imply  $g(\bar{T})=g(\bar{T}')$. Also, $g(T_n)=g(T_n')$ implies $g(\Phi^{m(n)}(T_n))=g(\Phi^{m(n)}(T_n'))$, so that 
$$\lim g(\Phi^{m(n)}(T_n'-x_n)) \stackrel{(7)}{=} g(\bar T') = g(\bar T) \stackrel{(6)}{=} \lim g(\Phi^{m(n)}(T_n)) = \lim g(\Phi^{m(n)}(T_n')).$$
Since, by (5),  $\lim g(\Phi^{m(n)}(T_n'-x_n)) = \lim g(\Phi^{m(n)}(T_n'))-x$ we find a contradiction to  
the freeness of the action of $x$ on $\mathbb{X}$.
\end{proof}

Remark:
Putnam proves in \cite{P} that if $g$ is any $u$-resolving factor map between Smale spaces (that is, $g$ is injective on unstable sets), then $g$ is finite-to-one. The above result is a corollary, under the assumption that $g$ is a semi-conjugacy with a suitably hyperbolic action on $\mathbb{X}$, as will be the case for the applications of Theorem \ref{thm-thm1} in this article.

The following is the extension of Cor.~\ref{finitely many pairs}.
 
\begin{corollary}\label{cor-finitely many pairs}
Assuming the hypotheses of Theorem \ref{thm-thm1}, up to translation, there are only finitely many
pairs of patches of the form $(B_0[T],B_0[T'])$ with $g(T)=g(T')$.
\end{corollary}

\begin{corollary}\label{cor-finite-to-one}
Assuming the hypotheses of Theorem \ref{thm-thm1} together with $g(T)=g(T')\implies g(\Phi^{-1}(T))=g(\Phi^{-1}(T'))$, $g$ is boundedly finite-to-one.
\end{corollary}

\begin{proof} Literally the same as for Thm.~\ref{finite-to-one}. 
\end{proof}

We now identify $\hat {\mathbb R}^n$ with the dual vector space ${\mathbb R^n}^*$ 
such that $b^*\in{\R^n}^*$ corresponds to the character $t\mapsto e^{2\pi\imath b^*(t)}$.
Then each endomorphism $\Lambda$ on $\mathbb R^n$ has a dual endomorphism which we denote $\Lambda^*$; w.r.t.\ an ONB it is the transpose of $\Lambda$. 
\begin{theorem}
Consider a Pisot family substitution of degree $d$ and multiplicity $J$ and with linear expansion $\Lambda$.
There exists a lattice $\Gamma$ of rank $dJ$ which is relatively dense in ${\R^n}^*$ such that
the group of eigenvalues of the $\R^n$-action is exactly
$$ \mathcal E = \dir (\Gamma,\Lambda^*). $$
\end{theorem}
\begin{proof}  Lee and Solomyak (\cite{LS}) show that for Pisot family substitutions 
there exist $J$ 
vectors $b_1^*,\cdots,b_J^*\in{\R^n}^*$ such that 
$\{{\Lambda^*}^m b^*_i:1\leq i\leq J, 0\leq m\leq d-1\}$ is a collection of eigenvalues  which is linear 
independent\footnote{Linear independence over $\mathbb Q$ is not stated explicitly in \cite{LS} but
can easily derived from what is written there.}
 over $\mathbb Q$ and spans ${\R^n}^*$. Let $\Gamma'$ be the group generated by these eigenvalues. It is a lattice of rank $dJ$. Since $\E$ is invariant under $\Lambda^*$ it contains the group $H$ 
generated by $\{{\Lambda^*}^m b^*_i:1\leq i\leq J, m\in\Z\} = \dir(\Gamma',\Lambda^*)$.
The strategy is to show that $i:H\hookrightarrow \mathcal E$ is a finite index inclusion of $H$ in $\E$. This is equivalent to saying that the map $\hat i:\hE \to \hat H$ is finite to one. This can be seen as follows.

Let $g:X\to \hat H$, $g=\hat i\circ\gmax$. This is a factor map as $g$ is surjective and $\R^n$-equivariant where the $\R^n$-action on $j\in \hat H$ is given by 
$( t\cdot j) ({\Lambda^*}^mb_i^*) = e^{2\pi \imath {\Lambda^*}^mb_i^*(t)} j ({\Lambda^*}^mb_i^*)$.
Since $H$ is relatively dense in $\hat G$,  
Lemma~\ref{lem-free} implies that this action is locally free.
Furthermore, ${\Lambda^*}^m\circ i = i \circ {\Lambda^*}^m$ so that
 ${\Lambda^*}^m\circ g = g \circ {\Phi}^m$. Hence $g$ satisfies the assumptions of Cor.~\ref{cor-finite-to-one}.
 It follows that $\hat\imath$ must be finite to one and $H$ a finite index subgroup of $\mathcal E$.

Now it follows that there are finitely many vectors $w_1,\cdots,w_k$ such that $\sum_i w_i + H = \E$.
It follows that for some $N$, $Nw_i\in H$ for all $i$ (otherwise the index would be infinite). By definition of the direct limit this means that for some $m$, $Nw_i\in {\Lambda^*}^{-m}\Gamma'$, and hence also $w_i\in \frac1N {\Lambda^*}^{-m}\Gamma'$. Let $\Gamma$ be the group generated by ${\Lambda^*}^{-m}\Gamma'$ and the $w_i$. Then we have inclusions $\Gamma'\subset\Gamma\subset \frac1N\Gamma'$ all of finite index. Hence 
$\Gamma$ has the same rank as $\Gamma'$. Note also that $\Gamma$ is invariant under $\Lambda^*$, as both $\E$ and $H$ are invariant. Thus we have $\Gamma\subset\E$ and $\E\subset \dir (\Gamma,\Lambda^*)$. Since $\E$ is invariant under $\Lambda^*$ this implies that the last inclusion is an equality. 
\end{proof}
\begin{cor}
The maximal equicontinuous factor is an inverse limit of $dJ$-tori, 
$$\hat{\mathcal E} = \inv (\mathbb T^{dJ},\widehat{\Lambda^*}), $$
$\mathbb T^{dJ} = \hat \Gamma$.
Its $\R^n$-action is 
free and $\Phimax$ is ergodic with respect to Haar measure.
\end{cor}
\begin{proof}
Only the last point needs a comment: By definition, a linear expansion has no eigenvalues on the unit circle, so the result follows from Lemma \ref{lem-ergodic}.\end{proof}
We note that if $\Phi$ admits a fixed point (which we can always achieve by going over to a power) $\Phimax=\widehat{\Lambda^*}$.  Here $\widehat{\Lambda^*}$ is the dual map to $\Lambda^*$ which is not to be confused with $\Lambda^{**} = \Lambda$ since dualisation is w.r.t.\ the group $\Gamma$ and not $\hat\R^n$.
In fact, by the Pisot family condition, $\widehat{ \Lambda^*}$ can be written  $\mbox{\rm diag}(A,\cdots,A)$ ($J$ copies) in some basis for some integer matrix $A$ whose characteristic polynomial is the minimal monic polynomial having the eigenvalues of $\Lambda$ as roots.
Note that if $\det A = \pm 1$ then $\E=\Gamma$ and $\hat{\mathcal E} = \mathbb T^{dJ}$.

\subsection{Further results}

The following is a generalisation of a result from \cite{BK} which is based on the definition of the coincidence rank.
\begin{lemma} \label{notprox} Consider a Meyer
  substitution tiling system, $x\in \Xmax$ and
  $T,T'\in \gmax^{-1}(x)$. If $T$ and $T'$ are not proximal then they
do not have a single tile in common:
$T\cap T'=\emptyset$. 
\end{lemma}
\begin{proof} Recall that $\delta_0$ is such that $cr =
  cr(x,\delta_0)$. If
$T,T'\in \gmax^{-1}(x)$ and $T$ and $T'$ are not proximal we thus have
$\inf_vd(T-v,T'-v)\geq\delta_0$ and $\delta_0$ does not depend on $T,T'$.
By Theorem~\ref{thm-Pisot-Meyer} and
Corollary~\ref{finitely many pairs} we can reformulate this as
$$\sup_v\{R:B_R[T-v] = B_R[T'-v])\}\leq R_0<+\infty$$ 
with $R_0$ not depending on $T,T'$. 
Since the substitution $\Phi$ preserves fibres of $\gmax$ and respects the proximality relation, we also have $\sup_v\{R:B_R[\Phi(T)-v] =
B_R[\Phi(T')-v])\}\leq R_0$. But $B_R[T-v] =
B_R[T'-v]$ implies that $B_{\lambda R}[\Phi(T)-v] =
B_{\lambda R}[\Phi(T')-v]$ for some $\lambda>1$. This is only possible
if $R_0=0$. 
\end{proof}
The coincidence rank thus counts the maximal number of tilings in a
fiber of $\gmax$ which are pairwise non-coincident in the sense that
they do not share a single tile. 

\section{Complete proximality}
\subsection{General results}
Consider a locally compact abelian group $G$ acting continuously on a compact metric space $X$.
Let $A\subset G$ be a subset which contains a translate of every compact subset of $G$.
We denote by $\mathcal A$ the collection of all such sets.

\begin{definition}
We say that  $x,y\in X$ are proximal in $A$, written $x \sim_{A,p}y$, if $$\inf_{t\in A} d(t\cdot x,t\cdot y) = 0$$ and $x,y$ are completely proximal, written 
$x \sim_{cp}y$, if $x,y$ are proximal in all $A\in\mathcal A$.
\end{definition}
For $\epsilon>0$ and $x,y\in X$ let $G_\epsilon(x,y) = \{ t\in G:d(t\cdot x,t\cdot y)\leq \epsilon\}$. Recall that a subset $B\subset G$ is called syndetic\footnote{For $G=\R^n$ the notions of  syndetic and relatively dense coincide.}  if there exists a compact subset $K\subset G$ such that $K+B=G$.
In particular, $G_\epsilon(x,y)$ is syndetic if and only if its complement $G_\epsilon(x,y)^c\notin \mathcal A$.
Note that $x\sim_{A,p} y$ whenever for all for all $\epsilon>0$, $G_\epsilon(x,y)\cap A\neq \emptyset$. From this
one sees easily:
\begin{lemma}
$x\sim_{cp} y$ if and only if for all $\epsilon>0$, $G_\epsilon(x,y)$ is syndetic.
\end{lemma}
\begin{proof}
If $G_\epsilon(x,y)$ is not syndetic then $G_\epsilon(x,y)^c\in \mathcal A$. Hence for $A=G_\epsilon(x,y)^c$ we can't have $x\sim_{A,p} y$. If $G_\epsilon(x,y)$ is syndetic then $G_\epsilon(x,y)^c\notin \mathcal A$ and hence there exists a compact $K$ such that $G_\epsilon(x,y)^c$ does not contain any of its translates. Hence no 
$A\in\mathcal A$ is contained in $G_\epsilon(x,y)^c$. But then all $A\cap G_\epsilon(x,y)\neq\emptyset$.
Since $\epsilon$ is arbritrary in that argument the statement follows.
\end{proof}
\begin{theorem} Suppose that $G$ is also $\sigma$-compact.
If proximality is closed, then it agrees with complete proximality.
\end{theorem}
\begin{proof} Clearly $x\sim_{cp} y$ implies $x\sim_p y$.
Suppose now that $x\sim_p y$ are two proximal elements of $X$ and let $A\in\A$. $G$ being $\sigma$-compact means that there exists a sequence $(K_n)_n$ of compact subsets such that 
$K_n\subset K_{n+1}$ and $\bigcup_n K_n = G$. 
Let $(t_n)_n$ be such that $t_n+K_n \subset A$. We may suppose (upon taking a subsequence) that
$\overline{x}=\lim t_n\cdot x$ and $\overline{y}=\lim t_n\cdot y$ exist. Since proximality is closed we have $\overline{x}\sim_p \overline{y}$ and hence there exists a sequence $(s_k)_k$ such that
$$ d(\lim_n s_k\cdot(t_n\cdot x),\lim_n s_k\cdot(t_n\cdot y)))<\frac1k.$$
For all $k$, there exists $n_k$ such that  $ d(s_k\cdot(t_{n_k}\cdot x),s_k\cdot(t_n\cdot x))<\frac1k$ for all $n\geq n_k$,
and similarily for $y$. We may also suppose that $s_k\in K_{n_k}$. Then $v_k=t_{n_k}+s_k\in A$ and
$d(v_k\cdot x,v_k\cdot y)<\frac3k$.
\end{proof}

\begin{theorem} Complete proximality is an equivalence relation.
\end{theorem} 
\begin{proof} Suppose that $x\sim_{cp}y$ and $y\sim_{cp}z$. Let $\epsilon>0$ be given. We will show that $G_\epsilon(x,z)$ is syndetic. 

$y\sim_{cp}z$ implies that there is a compact $K_1\subset G$ such that $\forall s\in G\exists k_1(s)\in K_1$:
$d((k_1(s)+s)\cdot y,(k_1(s)+s)\cdot z)\leq\frac\epsilon2$.

Since $K_1$ is compact and the group action, which we denote here by $\alpha$, continuous 
the family $\{\alpha_s:s\in K_1\}$ is uniformly equicontinuous. Hence there exists a $\delta>0$ such that $d(x_1,x_2)\leq \delta$ implies for all $s\in K_1$: $d(s\cdot x_1,s\cdot x_2)\leq \frac\epsilon2$.

 $x\sim_{cp}y$ implies that there is a compact $K_2\subset G$ such that $\forall t\in G\exists k_2\in K_2$:
$d((k_2+t)\cdot x,(k_2+t)\cdot y)\leq\delta$, and hence in particular, $\forall t\in G$ $\exists k_2\in K_2$ $\forall k_1\in K_1: d((k_1+k_2+t)\cdot x, (k_1+k_2+ t)\cdot y)\leq \frac\epsilon2$.

We claim that $K_1+K_2+G_\epsilon(x,z)=G$.
Indeed, given $t\in G$ there exists $k_2\in K_2$ such that $\forall k_1\in K_1: d((k_1+k_2+t)\cdot x,( k_1+k_2+ t)\cdot y)\leq \frac\epsilon2$. We choose $k_1=k_1(k_2+t)$ to obtain $d((k_2+k_1+t)\cdot y,(k_2+k_1+t)\cdot z)\leq \frac\epsilon2$ and hence $d((k_2+k_1+t)\cdot x,(k_2+k_1+t)\cdot z)\leq \epsilon$.
\end{proof}

Note that this gives a direct proof for Theorem~\ref{thm-Aus1} in the context of locally compact, $\sigma$ compact, abelian groups.

\begin{lemma}\label{lem-cp}
Let $x,y\in G$ be completely proximal and $\overline{x}=\lim t_n\cdot x$, $\overline{y}=\lim t_n\cdot y$ for some sequence $(t_n)_n\subset G$. Then $\overline{x}$ and $\overline{y}$ are completely proximal.
\end{lemma}
\begin{proof}
Let $\epsilon>0$ be given and $K$ compact such that $K+G_{\epsilon}(x,y)=G$. 
Pick any $s\in G$. Using uniform equicontinuity of the family $\{\alpha_s:s\in K\}$ we can find
an $n$ large enough so that $d((t_n+v)\cdot x,v\cdot\overline{x})\leq\epsilon$ and $d((t_n+v)\cdot y,v\cdot \overline{y})\leq\epsilon$ for all $v\in K$. There is $v\in K$ so that $d((t_n+v)\cdot x,(t_n+v)\cdot y)\leq\epsilon$ and hence
$d(v\cdot\overline{x}, v\cdot\overline{y})\leq 3\epsilon$.
\end{proof}
\subsection{Complete proximality for Meyer substitutions} In this section, $\Phi$ will be a substitution whose tiling space $\Omega$ has the Meyer property.
We show that complete proximality is a closed relation on $\Omega$. The group $G$ is $\R^n$ and we denote the action again by $T-v$ instead of $v\cdot T$. 

We say that a pair of patches $(P,P')$ occurs in the pair of tilings $(T,T')$ if there is $v\in\R^n$ so that $P-v\subset T$ and $P'-v\subset T'$. Let $\Rr$ be the relation on $\Omega$ defined by: $(S,S')\in \Rr$ if and only if for each $v\in\R^n$ and each $r>0$ there are $T,T'\in \Omega$ with $T\sim_{cp}T'$ so that $(B_r[S-v],B_r[S'-v])$ occurs in $(T,T')$.
\begin{lemma}\label {lem2} $\Rr$ is closed and $\Rr\subset \Qr$.
\end{lemma}
\begin{proof} We show first that $\Rr$ is closed: Suppose $(S_n)_n$ and $(S'_n)_n$ are two converging sequences, towards $S$ and $S'$, respectively, such that $(S_n,S'_n)\in \Rr$. Let $v$, $r$ be given. There exists 
$N$ such that for all
$n\geq N$ exists $\epsilon_n,\epsilon'_n$ such that
$$B_r[S-v-\epsilon_n] = B_r[S_n-v],\quad  B_r[S'-v-\epsilon'_n] = B_r[S'_n-v].$$
Furthermore, we may assume that the $\epsilon_n$ and $\epsilon'_n$ form sequences tending to $0$.
Since $(S_n,S'_n)\subset \Rr$ there exists $T,T'$ such that 
$ (B_r[S-v-\epsilon_n], B_r[S'-v-\epsilon'_n])$ occurs in $(T-\epsilon_n,T'-\epsilon'_n)$ and 
$T-\epsilon_n\sim_{cp} T'-\epsilon'_n$. By the Meyer property (Cor.~\ref{finitely many pairs}) we must have $\epsilon_n=\epsilon'_n$ for large enough $n$ and hence $T\sim_{cp} T'$. Furthermore
$ (B_r[S-v], B_r[S'-v])$ occurs in $(T,T')$ and hence $(S,S')\in \Rr$.

Let $(S,S')\in\Rr$.
For each $k\in\N$ there are $T_k,\sim_{cp}T'_k$ such that $(B_k[S],B_k[S'])$ occurs in $(T_k,T'_k)$. Then
$T_k\to S$ and $T'_k\to S'$. Since $T_k,\sim_{cp}T'_k$ implies $(T_k,T'_k)\in\Qr$ and $\Qr$ is closed $(S,S')\in\Qr$. 
\end{proof}

\begin{lemma}\label{lem1} Suppose that $(T,T')\in\Qr$, and suppose that the pair of finite patches $(P,P')$ occurs in $(T,T')$. Then there are $l\in\Z$ and $L\in\N$ so that $(P,P')$ occurs in $(\Phi^{-(kL+l)}(T),\Phi^{-(kL+l)}(T'))$ for all $k\in\N$.
\end{lemma}
\begin{proof}
Let $r>0$ be big enough so that $P\subset B_r[T]$ and $P'\subset B_r[T']$. Since $\Phi^{-1}$ preserves regional proximality we have
$(\Phi^{-k}(T),\Phi^{-k}(T'))\in\Qr$ for all $k$, and so we can conclude from Cor.~\ref{finitely many pairs} that
there are $k_i\to\infty$ so that 
$$(Q,Q'):=(B_r[\Phi^{-k_i}(T)],B_r[\Phi^{-k_i}(T')])$$ does not depend on $i$. 
Let $i$ be large enough so that $$B_r[T]\subset \Phi^{k_i}(
B_r[\Phi^{-k_i}(T)]),\quad B_r[T']\subset \Phi^{k_i}(
B_r[\Phi^{-k_i}(T')]).$$
Let $j>i$ be large enough so that  $$B_r[\Phi^{-k_i}(T)]\subset \Phi^{k_j-k_i}(
B_r[\Phi^{-k_j}(T)]),\quad B_r[\Phi^{-k_i}(T')]\subset \Phi^{k_j-k_i}(
B_r[\Phi^{-k_j}(T')]).$$
Let $L:=k_j-k_i$. There is $l'\in\{0,\ldots,L-1\}$ so that $k_s\equiv l'\mod{L}$ for infinitely many $s$.
Then $(Q,Q')$ occurs in $(\Phi^{-(kL+l')}(T),\Phi^{-(kL+l')}(T'))$ for all $k\in\N$ and $(P,P')$ occurs in 
 $(\Phi^{-(kL+l)}(T),\Phi^{-(kL+l)}(T'))$ for all $k\in\N$, where $l:=l'-k_i$.
\end{proof}

\begin{prop}  Complete proximality is closed for Meyer substitution tiling spaces.
\end{prop}
\begin{proof}
We will prove that $\Rr$ is the same as the complete proximality relation. That $T\sim_{cp}T'$ implies 
$(T,T')\in \Rr$ is immediate. So suppose that $(S,S')\in \Rr$. By Lemma~\ref{lem2},
$(S,S')\in\Qr$ and hence there is an increasing sequence of positive integers $k_i\to\infty$
so that  $$(P,P'):=(B_1[\Phi^{-k_i}(S)],B_1[\Phi^{-k_i}(S')])$$ 
does not depend on $i$. By recognizability (that is, since $\Phi$ is invertible), there is an $r>0$ so that if $$(B_r[T],B_r[T'])=(B_r[S],B_r[S'])$$ then $$(B_1[\Phi^{-k_1}(T)],B_1[\Phi^{-k_1}(T')])=(B_1[\Phi^{-k_1}(S)],B_1[\Phi^{-k_1}(S')]).$$
Since $(S,S')\in\Rr$ there are $T\sim_{cp}T'$ with $(B_r[T],B_r[T'])=(B_r[S],B_r[S'])$.  We 
apply Lemma \ref{lem1} to $(\Phi^{-k_1}(T),\Phi^{-k_1}(T'))$ to obtain $l\in\Z$ and $L\in\N$ so that $(P,P')$ occurs in
$$(\Phi^{-(kL+l)}(\Phi^{-k_1}(T)),\Phi^{-(kL+l)}(\Phi^{-k_1}(T')))$$ for all $k\in\N$. Let $l'\in\{0,\ldots,L-1\}$ be such that $k_i\equiv l'\mod{L}$ for infinitely many $i$: say $k_{i_j}=m_jL+l'$ with $m_j\in\N$ and $m_j\to\infty$. Let $\bar{T}:=\Phi^{l+k_1-l'}(T)$ and $\bar{T}':=\Phi^{l+k_1-l'}(T')$. Then $\bar{T}\sim_{cp}\bar{T}'$ and $(P,P')$ occurs in $(\Phi^{-k_{i_j}}(\bar{T}),\Phi^{-k_{i_j}}(\bar{T}'))$ for all $j$.
Say $B_1[\Phi^{-k_{i_j}}(S)]-v_j\in\Phi^{-k_{i_j}}(\bar{T}) ,$ and $B_1[\Phi^{-k_{i_j}}(S')]-v_j\in\Phi^{-k_{i_j}}(\bar{T}')$. Let $r_j\to\infty$ be such that $\Phi^{k_{i_j}}(B_1[\Phi^{-k_{i_j}}(S)])\supset
B_{r_j}[S]$ and $\Phi^{k_{i_j}}(B_1[\Phi^{-k_{i_j}}(S)])\supset B_{r_j}[S]$. Then $\bar{T}+\lambda^{k_{i_j}}v_j$ and  $\bar{T}'+\lambda^{k_{i_j}}v_j$ agree with $S$ and $S'$, resp., on $B_{r_j}(0)$, so that  $\bar{T}+\lambda^{k_{i_j}}v_j\to S$ and $\bar{T}'+\lambda^{k_{i_j}}v_j\to S'$ as $j\to \infty$,
and it follows from Lemma~\ref{lem-cp} that $S\sim_{cp}S'$.
\end{proof}
\begin{cor}\label{cor-Meyer}
For Meyer substitution tiling spaces proximality is closed if
and only if it coincides with complete proximality.
\end{cor}


\end{document}